\documentclass[10pt,twoside]{article}

\usepackage{latexsym}
\usepackage{amssymb}
\usepackage{array}
\usepackage{theorem}
\usepackage{pictex}

\newcommand{\CC}{\mathbb{C}}
\newcommand{\PP}{\mathbb{P}}

\newcommand{\TT}{\mathbb{T}}

\newcommand{\GG}{\mathbb{G}}

\newcommand{\ZZ}{\mathbb{Z}}
\renewcommand{\rho}{\varrho}
\renewcommand{\phi}{\varphi}
\renewcommand{\epsilon}{\varepsilon}
\newcommand{\Span}{\mathrm{span}\,}
\newcommand{\rank}{\mathrm{rank}\,}

\newcommand{\im}{\mathrm{Im}\,}
\newcommand{\sgn}{\mathrm{sgn}\,}

\newcommand{\id}{\mathrm{id}}
\renewcommand{\sec}{\mathrm{Sec}}
\newcommand{\aut}{\mathrm{Aut}}
\newcommand{\sym}{\mathrm{Sym}}
\newcommand{\antisym}{\mathrm{Antisym}}
\newcommand{\codim}{\mathrm{codim}\,}
\newcommand{\corank}{\mathrm{corank}\,}
\newcommand{\GL}{\mathrm{GL}}
\newcommand{\Sp}{\mathrm{Sp}}
\newcommand{\SL}{\mathrm{SL}}
\newcommand{\Pf}{\mathrm{Pf}}
\newcommand{\Tri}{\mathrm{Tri}}

\newcommand{\E}{\mathrm{E}}
\newcommand{\M}{\mathrm{M}}
\newcommand{\G}{\mathrm{G}}

\newcommand{\raiseghost}[3]{ % x y text
\raisebox{#2}[0cm][0cm]{\makebox[0cm]{\hspace{#1}{#3}}}}

\newtheorem{theorem}{Theorem}[section]

\newtheorem{remark}[theorem]{Remark}
\newtheorem{lemma}[theorem]{Lemma}
\newtheorem{corollary}[theorem]{Corollary}
\newtheorem{proposition}[theorem]{Proposition}

\begin{document} 

\title{The Automorphism Group of Linear Sections \\
of the Grassmannians $\GG(1,N)$}
\date{April 1999}
\author{J.\ Piontkowski and A.\ Van de Ven}
\maketitle

\begin{abstract} 
\noindent The Grassmannians of lines in projective $N$-space, $\GG(1,N)$, are 
embedded by
way of the Pl\"ucker embedding in the projective space 
$\PP(\bigwedge^2\CC^{N+1})$. Let $H^l$ be a general $l$-codimensional linear 
subspace in this projective space.

\noindent We examine the geometry of the linear sections $\GG(1,N)\cap H^l$ 
by studying
their automorphism groups and list those which are homogeneous or quasihomogeneous.  

\end{abstract}

\setcounter{section}{-1}
\section{Introduction}

Complete intersections in projective space have been studied extensively from
many points of view. A natural generalisation is the study of
 complete intersections
in Grassmannians. The first case that presents itself is the case of
intersections with linear spaces. Indeed, there is an extensive literature on
the simplest case, the Grassmannian of lines in the 3-space, where
intersections are known as linear complexes and congruences of lines. L.\ Roth
has studied the rationality of linear sections of Grassmannians of lines 
in general. 
If they are smooth and if the
dimension of the intersection is greater than half the dimension of the
Grassmannian, then they are rational. R.\ Donagi determined the cohomology
and the intermediate Jacobian of some linear sections of Grassmannians of 
lines.

\medskip
In this paper we study the linear sections from the point of automorphism
groups. Let $\GG(1,N)$ be the Grassmann variety 
of lines in projective $N$-space,
canonically embedded in $\PP({\textstyle\bigwedge}^2\CC^{N+1})$ and let $H^l$
be an $l$-codimensional linear subspace in this space. For general $H^l$ we
determine the automorphism groups for $\GG(1,N)\cap H$, $\GG(1,N)\cap H^2$,
$\GG(1,4)\cap H^3$, and $\GG(1,5)\cap H^3$. In the second case we find 
for example:

\setcounter{section}{3}
\setcounter{theorem}{4}
\begin{theorem}
For $N=2n-1\ge5$ the automorphism group of $\GG(1,N)\cap H^2$ has
$\SL(2,\CC)^n/\{1,-1\}$ as a normal subgroup and the quotient group is
isomorphic to the permutation group $\mathrm{S}(3)$ for $n=3$, to $\ZZ/2\ZZ
\times \ZZ/2\ZZ$ for $n=4$, and trivial otherwise.
\end{theorem}
\setcounter{section}{0}
\setcounter{theorem}{0}

\medskip
We believe that apart from trivial cases these are the only general linear
sections where automorphism groups of positive dimension appear. Extensive
computer checks seem to confirm this. 

\medskip
In particular we prove that the automorphism groups of
$\GG(1,2n)\cap H$, $\GG(1,4)\cap H^2$, $\GG(1,5)\cap H^2$, $\GG(1,6)\cap H^2$,
and $\GG(1,4)\cap H^3$
are quasihomogenous -- those of $\GG(1,2n-1)\cap H$ and $\GG(1,3)\cap H^2$ are
even homogenous -- whereas all others are not.

\medskip
As to our methods, in our proofs the rich geometry of the Grassmannian plays
a decisive role. Otherwise, we mainly use well known tools like multilinear
algebra, Lefschetz theorems, vanishing theorems etc.

\medskip
We are indebted to E.\ Opdam and A.\ Pasquale for useful remarks. The first
author also thanks the Stieltjes Institut of Leiden University for 
financial support.

\section{Preliminary}
The Grassmannian $\GG(1,N)$ of lines in $\PP_N$ is embedded by way of the
 Pl\"ucker embedding into $\PP(\bigwedge^2\CC^{N+1})$
\[ \begin{array}{ccc}
\GG(1,N) & \longrightarrow & \PP({\textstyle\bigwedge}^2\CC^{N+1})\\[1ex]
\Span\{ v,w \} &\longmapsto & \PP(v\wedge w).
\end{array} \]

We denote by $H^l$ an $l$-codimensional linear subspace of  
$\PP(\bigwedge^2\CC^{N+1})$. Roth \cite{R} 
examined the geometry of the general linear sections of the Grassmannians
 and found

\begin{theorem}
For  a general $H^l$ with $0 \le l \le 1/2 \dim \GG(1,N)=N-1 $ 
the intersection with  the 
Grassmannians, $\GG(1,N)\cap H^l$, is rational.
\end{theorem}

In this article we continue this study by describing the automorphism groups 
of these sections. As for the notation, given a subvariety $Y$ of a variety $X$
we define $\aut (Y,X)$ to be the automorphisms of $X$ that induce automorphisms
of $Y$, i.e.
\[
\aut(Y,X)=\{ \phi\in\aut(X)\mid \phi(Y)\subseteq Y \}.
\]

Recall that the automorphism group of the Grassmannian itself is computed in
two steps, see e.g. \cite[10.19]{H}. First one shows that all automorphisms 
are 
induced by 
automorphisms of  $\PP(\bigwedge^2\CC^{N+1})$, i.e.
\[
\aut (\GG(1,N))\cong 
\aut(\GG(1,N), \PP({\textstyle \bigwedge}^2\CC^{N+1})).
\]
Then one proves that for $N \not= 3$ the right hand side group is isomorphic 
to $\PP\GL(N+1,\CC)$ via
\[ \begin{array}{ccc}
\PP\GL(N+1,\CC) & \longrightarrow &
\aut(\GG(1,N),\PP({\textstyle\bigwedge}^2\CC^{N+1}))\\[1.5ex]
\PP(T) &\longmapsto & 
\left( \PP(\sum v_i\wedge w_i)\mapsto \PP(\sum Tv_i\wedge Tw_i)\right) .
\end{array} \]

For the linear sections of the Grassmannians we follow the same outline.
The first step is the following theorem; the second step will be done 
separately for the different cases in the next sections.

\begin{theorem} \label{autointersect} For a general linear subspace 
$H^l\subset \PP(\bigwedge^2\CC^{N+1})$ of codimension $l\le 2N-5$
\[
\aut(\GG(1,N)\cap H^l)=\aut(\GG(1,N),\PP({\textstyle \bigwedge}^2\CC^{N+1}))
\cap \aut(H^l,\PP({\textstyle \bigwedge}^2\CC^{N+1})).
\]
\end{theorem}
\emph{Proof.} We will abbreviate $\GG(1,N)$ by $\GG$. The ``$\supseteq$'' 
inclusion is trivial. For the other one we prove first that all automorphisms 
of $\GG\cap H^l$ are induced by automorphisms of $\GG$. This follows 
immediately, once we show that all divisors of $\GG\cap H^l$ are induced by 
divisors of $ \PP(\bigwedge^2\CC^{N+1})$, i.e.
\[ \mathrm{Pic}(\GG\cap H^l)=\mathrm{Pic}(\GG)=\ZZ\cdot H. \]

To see this, note that by the Lefschetz hyperplane section theorem
\[
\ZZ\cdot H=\mathrm{H}^2(\GG,\ZZ)=\mathrm{H}^2(\GG\cap H,\ZZ)=\ldots=\mathrm{H}^2(\GG\cap H^l,\ZZ)
\]
for $0\le l \le 2N-5$. From the exponential sequence
\[
0\rightarrow \ZZ_{\GG\cap H^l}\rightarrow {\cal O}_{\GG\cap H^l}
\rightarrow  {\cal O}_{\GG\cap H^l}^* \rightarrow 0
\]
we get as a part of the associated long exact sequence
\[
\ldots \rightarrow \mathrm{H}^1(\GG\cap H^l,{\cal O}) \rightarrow 
 \mathrm{H}^1(\GG\cap H^l,{\cal O}^*) \rightarrow  \mathrm{H}^2(\GG\cap H^l,\ZZ)=\ZZ\cdot H
\rightarrow 0
\]
and therefore
\[
\mathrm{Pic}(\GG\cap H^l)=\mathrm{H}^1(\GG\cap H^l,{\cal O}^*)=\ZZ\cdot H
\]
as soon as we know that $\mathrm{H}^1(\GG\cap H^l,{\cal O})=0$.

\medskip
This is well known for $l=0$.
For $l\ge 1$ we look at the restriction sequence
\[
0\rightarrow {\cal O}_{\GG\cap H^{l-1}}(-H) \rightarrow 
{\cal O}_{\GG\cap H^{l-1}} \rightarrow
{\cal O}_{\GG\cap H^l} \rightarrow 0
\]
and take its associated long exact sequence
\[ 
\begin{array}{l}
\displaystyle 
\ldots \rightarrow \mathrm{H}^1(\GG\cap H^{l-1},{\cal O}(-H)) \rightarrow 
 \mathrm{H}^1(\GG\cap H^{l-1},{\cal O}) \rightarrow  \mathrm{H}^1(\GG\cap H^l,{\cal O})
\rightarrow \\[1.5ex]
%\multicolumn{1}{r}{
\displaystyle
\phantom{\ldots}
\rightarrow\mathrm{H}^2(\GG\cap H^{l-1},{\cal O}(-H)) \rightarrow \ldots
%}
\end{array}
\]
The right  and left cohomology groups are trivial for $l\le 2N-4$ by Kodaira's 
vanishing theorem, so 
\[
0=\mathrm{H}^1(\GG,{\cal O})=\mathrm{H}^1(\GG\cap H,{\cal O})=\ldots=
\mathrm{H}^1(\GG\cap H^l,{\cal O}).
\]

\medskip
Now that we know $\aut(\GG\cap H^l)\subseteq \aut(\PP(\bigwedge^2\CC^{N+1}))$
it remains to show that these projective transformations of 
$\PP(\bigwedge^2\CC^{N+1})$ fix $H^l$. This will follow if we 
prove  that $\GG\cap H^l$ spans $H^l$, i.e.
\[
\mathrm{h}^0(\GG\cap H^l,{\cal O}(H))=\dim {\textstyle \bigwedge^2 }
\CC^{N+1}-l.
\]

This is known for $l=0$. For $l\ge 1$ we take the long exact sequence 
associated to the restriction sequence tensored by ${\cal O}(H)$
\[
\begin{array}{l}
\displaystyle 
0\rightarrow \mathrm{H}^0(\GG\cap H^{l-1},{\cal O})=\CC \rightarrow
\mathrm{H}^0(\GG\cap \mathrm{H}^{l-1},{\cal O}(H))\rightarrow
\mathrm{H}^0(\GG\cap \mathrm{H}^l,{\cal O}(H))\rightarrow \\[1.5ex]
\displaystyle
\phantom{0}
\rightarrow \mathrm{H}^1(\GG\cap \mathrm{H}^{l-1},{\cal O})=0.
\end{array}
\]
Looking at the dimensions we get 
\[
\mathrm{h}^0(\GG\cap H^l,{\cal O}(H))=
\mathrm{h}^0(\GG\cap H^{l-1},{\cal O}(H))-1,
\]
and the claim follows by induction. 
\hfill $\Box$

\bigskip

It is tempting to assume that the groups 
$\aut(\GG(1,N),\PP(\bigwedge^2\CC^{N+1}))$ and 
$\aut(H^l,\PP(\bigwedge^2\CC^{N+1}))$ in $\aut(\PP(\bigwedge^2\CC^{N+1}))$ 
intersect transversally.
Then the dimension of $\aut(\GG(1,N)\cap H^l)$ could be computed as
\[ \begin{array}{c@{\;=\;}l}
\dim \aut(\GG(1,N)\cap H^l) & \dim \aut (\GG(1,N))-
\codim \aut(H^l,\PP(\bigwedge^2\CC^{N+1}))\\[1ex]
 &(N+1)^2-1-l\left( { N+1 \choose 2}-l\right).
\end{array} \]
And we would find the following non-finite groups:
\[\begin{array}{l@{\;=\;}l}
\dim \aut(\GG(1,N)\cap H) & (N^2+3N+2)/2 \\[1ex]
\dim \aut(\GG(1,N)\cap H^2)& N+4 \\[1ex]
\dim \aut(\GG(1,4)\cap H^3)&3.
\end{array}\]

Unfortunately, the intersection is not always transversal.  Our computation 
of the automorphism groups will show the following dimensions for $N\ge 4$:
\[\begin{array}{l@{\;=\;}l}
\dim \aut(\GG(1,N)\cap H) & (N^2+3N+2)/2 \\[1ex]
\dim \aut(\GG(1,N)\cap H^2)& \left\{ \begin{array}{ll} 
                                      N+4& \mbox{for\ }N\mbox{\ even}\\
                                      3(N+1)/2 &  \mbox{for\ }N\mbox{\ odd}
                             \end{array}\right. \\[2.5ex]
\dim \aut(\GG(1,4)\cap H^3)&3\\[1ex]
\dim \aut(\GG(1,5)\cap H^3)&1.
\end{array}\]

We conjecture that these are the only non-finite groups. For $N+2 \le l $ the 
canonical bundle $K={\cal O}(-N-1+l)$ is positive on $\GG(1,N)\cap H^l$, and  
this conjecture can be proved by Serre's duality theorem and Kodaira's 
vanishing theorem:
\[ 
\dim \aut(\GG \cap H^l) = h^0(\GG \cap H^l, \Theta)=
h^{2N-2-l}(\GG \cap H^l, K \Omega^1)=0
\]
A proof for the remaining cases $3\le l \le N+1$ seems difficult. 
By computer computations we verified the conjecture for $N\le 10$ and all $l$.

\bigskip
With this theorem our task of determining the automorphisms of 
$\GG(1,N)\cap H^l$ has been immensely simplified. All we need to do is to find 
the projective transformations of  
$\aut(\GG(1,N),\PP(\bigwedge^2\CC^{N+1}))=\PP\GL(N+1,\CC)$ such 
that their induced action on $\PP(\bigwedge^2\CC^{N+1})^*$ preserves $H^l$.
To express this in algebraic terms we identify $(\bigwedge^2\CC^{N+1})^*$
with $\bigwedge^2(\CC^{N+1})^*$. If a particular basis of $\CC^{N+1}$ is
chosen,  $\bigwedge^2(\CC^{N+1})^*$ as antisymmetric forms on $\CC^{N+1}$
can also be identified with the antisymmetric matrices of size $N+1$. 
In concrete terms, if $(e_0,\ldots,e_N)$ is a basis of $\CC^{N+1}$ and
$E_{ij}\in \M(N+1,\CC)$ the matrix, which has a $1$ in the position $(i,j)$
but is otherwise zero, then 
\[
\begin{array}{ccc}
\left( {\textstyle \bigwedge^2}\CC^{N+1} \right)^* & \longrightarrow &
\antisym(N+1,\CC)\\[2ex]
\sum_{i,j} \lambda_{ij}(e_i\wedge e_j)^* &\longmapsto & 
{\textstyle \frac{1}{2}} \sum_{i,j} \lambda_{ij} (E_{ij}-E_{ji}).
\end{array}
\]

In these terms a line $l=p\wedge q\in \GG(1,N)$ is in the hyperplane 
$H\in \PP(\bigwedge^2\CC^{N+1})$ iff for a corresponding antisymmetric matrix 
$A\in \antisym(N+1,\CC)$ with $\PP(A)=H$ we have ${}^t\!p A q=0$.

\medskip
Further, the action of $\PP\GL(N+1,\CC)$ on $\PP(\bigwedge^2\CC^{N+1})$,
which was given for $\PP(T)\in\PP\GL(N+1,\CC)$ by 
\[
\begin{array}{ccc}
\PP\left( {\textstyle \bigwedge^2}\CC^{N+1} \right) & \longrightarrow &
\PP\left( {\textstyle \bigwedge^2}\CC^{N+1} \right)  \\[2ex]
\PP ( \sum v_i\wedge w_i )  &\longmapsto & 
\PP ( \sum Tv_i\wedge Tw_i ),
\end{array}
\]
induces the following action on the dual space
\[
\begin{array}{ccc}
\PP( \antisym(N+1,\CC) ) & \longrightarrow &
\PP( \antisym(N+1,\CC) ) \\[1ex]
\PP(A)  &\longmapsto & 
\PP ({}^tT^{-1} A T^{-1} ).
\end{array}
\]
Hence an $l$-codimensional linear subspace $H^l\subseteq\PP( \bigwedge^2\CC^{N+1})$ which is 
dually given
by $\PP(\Span\{ A_1,\ldots,A_l \})$ is preserved under $T$ iff every hyperplane
containing $H^l$ is mapped to another hyperplane containing $H^l$, i.e.
\[
\begin{array}{ll}
 {}^tT^{-1}\left(\sum \lambda_i A_i\right)T^{-1} \in
\Span \{ A_1,\ldots,A_l \} & 
 \mbox{ for\ all\ } \lambda_i\in\CC\\[2ex]
\Longleftrightarrow {}^tT^{-1}A_iT^{-1} \in
\Span \{ A_1,\ldots,A_l \} & \mbox{ for\ } i=1\ldots l.
\end{array}
\]
We conclude

\medskip
\begin{corollary}
For $N\ge 4$, $0\le l \le 2N-5$ and a general 
$H^l\subset \PP( \bigwedge^2\CC^{N+1})$ given by 
$\PP(\Span\{ A_1,\ldots,A_l \}) \subset  \PP(\antisym(N+1,\CC))$ the 
automorphism group of $\,\GG(1,N)\cap H^l$ is
\[
%\aut(\GG(1,N)\cap H^l)=
\{\PP(T)\in\PP\GL(N+1,\CC) \mid {}^t T^{-1} A_i T^{-1} 
\in \Span \{ A_1,\dots,A_l\} 
\, \forall i \}.
\]
\end{corollary}

%\emph{Proof.} Since we compute a group, it does not matter if we consider
%$T$ or $T^{-1}$ and the statement follows from what was said above. 
%\hfill $\Box$

\bigskip
In the following sections we will compute the automorphism groups using
this corollary. In the course of 
the computations we will use geometric arguments
for which it is essential to know if a hyperplane 
$H\subset \PP(\bigwedge^2\CC^{N+1})$ is tangent to $\GG(1,N)$ or not.
We recall the basic facts together with their short proofs.

\medskip
\begin{proposition} \label{schubert}
For any line $l_0\in\GG(1,N)$ the Schubert cycle
\[
\sigma:=\{ l\in \GG(1,N)\mid l\cap l_0\not=\emptyset \} \subseteq \GG(1,N)
\]
lies inside the tangent space 
$\TT_{l_0}\GG(1,N)\subseteq\PP(\bigwedge^2\CC^{N+1}) $
and spans it.
\end{proposition}
\emph{Proof.} Let $l\in\sigma$, $p\in l\cap l_0$, $q\in l_0\setminus \{p\}$ 
and 
$r\in l\setminus \{ p \}$ then 
\[
\begin{array}{ccc}
\CC & \longrightarrow & \GG(1,N) \\[1ex]
\lambda &\longmapsto & p\wedge (q+\lambda r)
\end{array}
\]
is a line in $\sigma\subset\GG(1,N) $ through $l_0$ and $l$. Therefore it is 
contained in the tangent space $\TT_{l_0}\GG(1,N)$, in particular 
$l\in\TT_{l_0}\GG(1,N)$.

\smallskip
We choose a basis $(e_0,\ldots,e_N)$ of $\CC^{N+1}$ such that 
$l_0=\PP(e_0\wedge e_1)$. The $2N-1$ points $\PP(e_0\wedge e_1)$,
$\PP(e_0\wedge e_i)$, $\PP(e_1\wedge e_i)$ for $i=2 \ldots N$ lie in
$\sigma \subset \TT_{l_0}\GG(1,N)$ and are projectively independent, 
hence they span  $\TT_{l_0}\GG(1,N)$.
\hfill$\Box$

\bigskip
\begin{corollary} \label{hyperplanetangency}
Let $H=\PP(A)\in\PP(\bigwedge^2\CC^{N+1})^*$ be a hyperplane and 
$l_0 \in \GG(1,N)$ a line then
\[
\TT_{l_0}\GG(1,N)\subseteq H \Longleftrightarrow l_0\subseteq \ker A.
\]
\end{corollary}
\emph{Proof.} By the Proposition $\TT_{l_0}\GG(1,N)\subseteq H $ is equivalent
to $\sigma \subseteq H$. 
If we use the same basis of $\CC^{N+1}$ as in the proof
of the proposition, this means that 
\[ \begin{array}{l}
\PP((\lambda e_0+\mu e_1)\wedge v)\in H \quad \mbox{for\ all\ }
(\lambda\colon\mu)\in \PP_1, \ v\in\CC^{N+1}\\[1ex]
\Longleftrightarrow {}^t\!(\lambda e_0+\mu e_1)  A  v =0  \quad \mbox{for\ all\ }
(\lambda\colon\mu)\in \PP_1, \ v\in\CC^{N+1}\\[1ex]
\Longleftrightarrow l_0 \subseteq \ker A.
\end{array}
\]

\vspace{-4.1ex}

\ \hfill$\Box$

\medskip
\begin{corollary}\label{dualG}
The dual variety $\GG(1,N)^*\subset \PP(\bigwedge^2\CC^{N+1})^*$ of
the Grassmannian variety $\GG(1,N)$ consists of matrices of corank $\ge 2$ for
$N$ odd resp. corank $\ge 3$ for $N$ even.

\medskip
For $N$ odd it is an irreducible hypersurface of degree $(N+1)/2$;
for $N$ even it is a 3-codimensional subvariety. 
\end{corollary}
\emph{Proof.} By the last corollary $H=\PP(A)\in \PP(\bigwedge^2\CC^{N+1})^*$ 
is
tangential to $\GG(1,N)$ iff $\corank A\ge 2$.
Recall that an antisymmetric matrix has even rank. So, for  $N$ odd the matrix 
$A\in \antisym(N+1,\CC)$ has corank $\ge 2$ iff $\det A=0$. But again since $A$
is antisymmetric, $\det A$ is the square of the irreducible Pfaffian polynomial
$\Pf A$ \cite[5.2]{B}, which therefore defines $\GG(1,N)^*$.

\medskip
For $N$ even $\corank A \ge 2$ is equivalent to $\corank A \ge 3$. We 
compute the dimension of $\GG(1,N)^*$ following Mumford \cite{M} and find
\[
\begin{array}{l}
\begin{array}{c@{\;=\;}l}
\dim\left( \begin{array}{l} \mbox{space\ of\ }A \mbox{\ with} \\ 
           \dim \ker A =3 \end{array} \right)
 & \dim \G(3,N+1) + \dim \bigwedge^2 \CC^{N+1}/\CC^3\\[-0.4ex]
 & 3(N-2) +(N-2)(N-3)/2\\[1ex]
 & (N^2+N-6)/2
\end{array}\\[7ex]
\Longrightarrow \codim \GG(1,N)^*=(N+1)N/2-(N^2+N-6)/2=3.
\end{array}
\]

\vspace{-4.3ex}

\ \hfill $\Box$

\section{$\GG(1,2n-1)\cap H$}\label{g12n-1hsection}

Let the hyperplane $H\in \PP(\bigwedge^2\CC^{2n})$ be given by an element 
$A \in (\bigwedge^2\CC^{2n})^*$, which we identify with its corresponding
antisymmetric matrix. If $H$ is general, $A$ will be a matrix of full rank.
This may be taken as the definition of a general $H$. We will assume from now 
on that $H$ is general. 

\medskip
The line system $\GG(1,2n-1)\cap H$ in $\PP_{2n-1}$ does not lead to obvious 
special points in the $\PP_{2n-1}$. Through every point  $p\in \PP_{2n-1}$
passes a $\PP_{2n-3}$ of lines, namely
\[
p\wedge q \in \GG(1,2n-1) \ \mbox{with}  \ q\in \ker {}^t\! pA.
\]

For $n\ge 3$ we can compute the automorphism group of $\GG(1,2n-1)\cap H$
with the help of Theorem \ref{autointersect} and its Corollary. 
It consists of elements
$\PP(T)\in \PP \GL(2n,\CC)=\aut( \GG(1,2n-1))$ such that $\PP(T)$ as an
element of $\PP\GL(\bigwedge^2\CC^{2n})$ preserves $H$, i.e.
\[ 
{}^t T^{-1} A T^{-1} = \lambda A \ \ 
\mbox{\ for\ suitable\ } \lambda \in \CC^*.
\]
We may choose coordinates on $\PP_{2n-1}$ such that 
\[ 
A=\left(
\begin{array}{cc}
0 & -\E_n \\
\E_n & 0
\end{array}
\right).
\]
Then by definition
\[ 
\Sp(2n,\CC)=\{ T\in \GL(2n,\CC) \mid {}^t T^{-1}A T^{-1}=A \}, 
\] 
and we have an isomorphism
\[
\begin{array}{c@{\;}c@{\;}c}
\{ T\in \GL(2n,\CC) \mid \exists\lambda_T\in \CC^* : 
{}^t T^{-1} A T^{-1}=\lambda_T A\}/\CC^* &
\longrightarrow & \Sp(2n,\CC)/\{1,-1\} \\[1.7ex]
\CC^*\cdot T & \longmapsto & \pm \frac{1}{\sqrt{\lambda_T}} T.
\end{array}
\]
Therefore we see

\begin{proposition} The automorphism group of $\GG(1,2n-1)\cap H$
for a general $H\subset \PP(\bigwedge^2\CC^{N+1})$ is  
$\Sp(2n,\CC)/\{1,-1\}$. Its action on  $\GG(1,2n-1)\cap H$ is homogeneous.
\end{proposition}
\emph{Proof.} The missing case of $\GG(1,3)\cap H$ can be found in 
\cite[p. 278]{FH}. The transitivity of the action follows from 
Witt's theorem \cite[12.31]{Br}.
\hfill $\Box$

\section{$\GG(1,2n-1)\cap H^2$}

A 2-codimensional linear subspace $L=H^2$ of $\PP(\bigwedge^2\CC^{2n})$
can be thought of as the pencil of hyperplanes containing it. So it gives
a line $L^*=\PP(\lambda A-\mu B)\subset \PP(\bigwedge^2\CC^{2n})^*$.
We identify again $(\bigwedge^2\CC^{2n})^*$ with the antisymmetric matrices
of size $2n$. 
The line $L^*$ intersects the dual Grassmannian $\GG(1,2n-1)^*$,
which consists of antisymmetric matrices of rank $\le 2n-2$ 
and is a hypersurface of degree $n$ by Corollary \ref{dualG},
in at most $n$ points.
For the moment a line $L^*$, and hence $L$, will be called general 
if it has $n$ points of 
intersection, $H_i=\PP(\lambda_i A-\mu_i B)\in L^*$, $i=1\ldots n$, with
the dual Grassmannian.
These hyperplanes 
$H_i$ are tangent to the Grassmannian $\GG(1,2n-1)$ at the points 
$l_i:= \ker (\lambda_i A-\mu_i B)\in\GG(1,2n-1)$ by Corollary \ref{hyperplanetangency}. Therefore  we get $n$ 
exceptional lines $l_1,\ldots,l_n$ in $\PP_{2n-1}$.

\medskip
The intersection of the Grassmannian with its tangent hyperplane $H_i$ 
contains all lines that intersect $l_i$, because these lines are already 
contained in the intersection $\GG(1,2n-1)\cap \TT_{l_i}\GG(1,2n-1)$
by Proposition \ref{schubert}.

\medskip
So, any line through a point $p\in l_i$ will be in the subspace 
$L\subset  \PP(\bigwedge^2\CC^{2n})$ as soon as it is contained
in any other hyperplane
$H\in L^*\setminus\{H_i\}$. This gives one linear restriction to lines 
through $p$, so that there is at least a $\PP_{2n-3}$ of lines through  the
points of the lines $l_i$. In contrast, through a general point of 
$\PP_{2n-1}\setminus \bigcup l_i$ there is only a $\PP_{2n-4}$ of lines.
In fact, we have 

\begin{proposition}\label{linecharacter}
The points of the lines $l_1,\ldots,l_n$ are characterized by the property
that through each of them passes a $\PP_{2n-3}$ of lines, i.e.
\[
\left\{ p\in \PP_{2n-1}\;
\begin{array}{|l}
\mbox{through\ } p \mbox{\ passes\  a\ }\PP_{2n-3}\\
\mbox{of\ lines\ of\ } \GG(1,2n-1)\cap L
\end{array}\right\}
=\bigcup l_i.\]
Furthermore, the lines $l_1,\ldots,l_n$ span the whole $\PP_{2n-1}$.
\end{proposition}

This may easily be seen if we write the pencil of hyperplanes $L^*$ in its 
normal form.

\begin{proposition}[Donagi\cite{D}]\label{normalpencil}
Given a pencil of hyperplanes $L^*=\PP(\lambda A-\mu B)\subset \PP(\bigwedge^2\CC^{2n})^*$ such
that the line $L^*$ intersects the Pfaffian hypersurface in $n$ different points.
Then there is a basis of $\CC^{2n}$ such that
\[
A=\left( \begin{array}{ccc}
J &  &  \raiseghost{-2ex}{-1ex}{{\emph{\LARGE 0}}}\\
  &\ddots & \\
\raiseghost{2ex}{0ex}{{\emph{\LARGE 0}}} & &J
\end{array} \right)
\ \ \mbox{and} \ \
B=\left( \begin{array}{ccc}
\lambda_1 J &  & \raiseghost{-2ex}{-1ex}{{\emph{\LARGE 0}}}\\
  &\ddots & \\
\raiseghost{2ex}{0ex}{{\emph{\LARGE 0}}} & &\lambda_n J
\end{array} \right)
\ \ \mbox{with} \
J=\left( \begin{array}{cc}
0 &-1  \\
1 & 0
\end{array} \right)\!.
\]
The points 
$(\lambda_1\,\colon 1),\ldots,(\lambda_n\,\colon 1)\in \PP_1\cong L^*$ 
are unique up to a projective transformation of $\PP_1$.
\end{proposition}

\noindent\emph{Proof of Proposition \ref{linecharacter}.} 
The hyperplane $H_i=\PP(\lambda_i A-\mu_i B)$ has, 
written as an antisymmetric
matrix, the kernel $l_i=\Span\{e_{2i-1},e_{2i}\}$ which means it is tangent to 
$\GG(1,2n-1)$ at $l_i$. All lines of $\GG(1,2n-1)\cap L$ through the point $p$
are given by $p\wedge q$ with ${}^t\!pAq={}^t\!pBq=0$. In order to have a 
$\PP_{2n-2}$ of lines through $p$,
the linear forms ${}^t\!pA$ and ${}^t\!pB$ must be
linear dependent, i.e. there are $\lambda,\mu\in\CC$ with
\[ 
0=\lambda {}^t\!pA -\mu {}^t\!pB={}^t\!p (\lambda A -\mu B).
\]
Therefore $p$ is in the kernel of a matrix of the pencil, but these kernels are
the lines $l_1,\ldots,l_n$, so $p$ is contained in one of them.
\hfill$\Box$

\bigskip
Knowing the  exceptional lines $l_1,\ldots,l_n$, one can immediately give some lines
which are in the line system.

\begin{proposition}
Any line in $\PP_{2n-1}$ which intersects two exceptional lines is an element
of the line system $\GG(1,2n-1)\cap L$.

\medskip
\noindent The exceptional lines themselves are not in the line system.
\end{proposition}
\emph{Proof.} If a line $l$ intersects $l_i$ and $l_j$, it lies 
-- as a point of the Grassmannian $\GG(1,2n-1)$ -- 
in $H_i$ and $H_j$, hence in $L=H_i\cap H_j$.

\medskip
Assume that the exceptional line $l_1$ is an element of $\GG(1,2n-1)$. 
By Proposition \ref{linecharacter} the lines through a point $p \in l_1$
sweep out a hyperplane. This hyperplane contains the line $l_1$ by assumption
and the other exceptional lines $l_2,\ldots,l_n$ by  the first part of this
proposition. But this contradicts the second statement of 
Proposition \ref{linecharacter}.
\hfill$\Box$

\bigskip
\begin{remark} From Proposition \ref{normalpencil} we also see that 
any position
of the $n$ points of the line $L^*$ is possible. In particular, 
we may call a 
line general if the position of the points is general 
in the sense needed below. 
\end{remark}

\bigskip
Using this geometric description we can determine the automorphisms of 
$\GG(1,2n-1)\cap L$. For the moment we restrict ourselves to $n\ge 3$ in
order to be able to use Theorem \ref{autointersect}. 
By this theorem  and its Corollary we can view an automorphism of 
$\GG(1,2n-1)\cap L$ 
as an element $\PP(T)$ of $\PP\GL(2n,\CC)$. 
To make the notation simpler, we will 
write only $T$ for $\PP(T)$ if no confusion can result. Since the points of the
exceptional lines are characterized by the property of Proposition 
\ref{linecharacter}, $T$ must map the union of the lines 
$l_i\subset \PP_{2n-1}$ onto itself. Permutations of the lines may occur, 
but -- as we will presently see --  not all permutations are possible.

\medskip
If we view the automorphism $T$ as an element of
$\aut (L,\PP({\textstyle \bigwedge\nolimits^2}\CC^{2n}))$,
it interchanges the hyperplanes containing $L$, i.e.\ it induces a projective
transformation of the line $L^*\subset \PP(\bigwedge\nolimits^2\CC^{2n})$.
Naturally,  the transformation of $L^*$ must preserve the union of points
of intersection of $L^*$ with the dual Grassmannian, which determine the 
lines $l_i$. Now, if a transformation of $\PP_{2n-1}$ permutes the lines $l_i$,
then the induced transformation of $L^*$ must permute the corresponding points
of $L^*$ in the same way.

\medskip
Since not every permutation of four or more points on a line can be induced by
a projective transformation, not all permutations are possible. In fact, if
the points are in general position, we get the following subgroups of the 
permutation groups:
\[
\begin{array}{c|c}
n & \mbox{subgroup\ of\ } \mathrm{S}(n)\\[0.5ex] \hline
3 & \mathrm{S}(3)\\[0.5ex]
4 & \{ (1\,2\,3\,4),(2\,1\,4\,3),(3\,4\,1\,2),(4\,3\,2\,1)\}\cong 
\ZZ/2\ZZ\times\ZZ/2\ZZ \\[0.5ex]
\ge 5& \{\id\}
\end{array}
\]

On the other hand, any permutation $\sigma\in\mathrm{S}(n)$ of the points on $L^*$
that is induced by a projective transformation $\phi$ of $L^*$ can be induced
by an automorphism of $\GG(1,2n-1)\cap L$. To see this, let us write $L^*$ in 
its normal form and define $T\in \GL(2n,\CC)$ as
\[
T(e_{2i}):=e_{2\sigma(i)} \quad \mbox{and} \quad
T(e_{2i-1}):=e_{2\sigma(i)-1}.
\]
This transformation permutes the lines in the prescribed way, and as an 
automorphism of $\PP(\bigwedge^2\CC^{2n})$ it fixes $L$ since the transformed
line $L^*$ is
\[
\begin{array}{l}
{}^tT^{-1}
\left( \lambda 
\left( \begin{array}{ccc}
J &  & \raiseghost{-2ex}{-1ex}{{\LARGE 0}}\\
  &\ddots & \\
\raiseghost{2ex}{0ex}{{\LARGE 0}} & & J
\end{array} \right)
-\mu
\left( \begin{array}{ccc}
\lambda_1 J &  &   \raiseghost{-2ex}{-1ex}{{\LARGE 0}} \\
  &\ddots & \\
\raiseghost{2ex}{0ex}{{\LARGE 0}} & &\lambda_n J
\end{array} \right)
\right)T^{-1} \\[6ex]
=   \lambda 
\left( \begin{array}{ccc}
J &  &  \raiseghost{-2ex}{-1ex}{{\LARGE 0}}\\
  &\ddots & \\
\raiseghost{2ex}{0ex}{{\LARGE 0}} & & J
\end{array} \right)
-\mu
\left( \begin{array}{ccc}
\lambda_{\sigma^{-1}(1)} J &  &  \raiseghost{-2ex}{-1ex}{{\LARGE 0}}\\
  &\ddots & \\
\raiseghost{2ex}{0ex}{{\LARGE 0}} & &\lambda_{\sigma^{-1}(n)} J
\end{array} \right).
\end{array}
\]

Changing the parametrisation of the line by $\phi$ we get back the old 
parametrisation of the line $L^*$ by the definition of $\phi$. 
So this $T$ is an
automorphism of $\GG(1,2n-1)\cap L$ that induces the permutation of lines we
started with.

\bigskip
Now we can restrict our attention to transformations that do not permute the 
lines since we can obtain every permutation by composing with one of the 
transformations from above. A transformation leaving all the lines individually fixed has the form
\[ 
T=
\left( \begin{array}{ccc}
t_1 &  &  \raiseghost{-2ex}{-1ex}{{\LARGE 0}}\\
  &\ddots & \\
\raiseghost{2ex}{0ex}{{\LARGE 0}} & &t_n
\end{array} \right)
\ \ \mbox{with\ }t_1,\ldots,t_n\in\GL(2,\CC).
\]

This $T$ will fix the line system $\GG(1,2n-1)\cap L$ in $\PP_{2n-1}$
iff it preserves $L^*$, i.e.\ for all $\lambda,\mu\in\CC$ there exists 
$\alpha,\beta\in\CC$ such that
\[ 
{}^tT^{-1}(\lambda A-\mu B)T^{-1}=\alpha A+\beta B.
\]
It is sufficient to check this for $(\lambda,\mu)=(1,0)$ and $(0,-1)$. Since
\[
\begin{array}{c@{\,}c@{\,}c}
{}^tT^{-1}AT^{-1} & = &
\left( \begin{array}{ccc}
\det t_1^{-1} J &  &  \raiseghost{-2ex}{-1ex}{{\LARGE 0}}\\
  &\ddots & \\
\raiseghost{2ex}{0ex}{{\LARGE 0}} & &\det t_n^{-1} J
\end{array} \right)\\[5.5ex]
{}^tT^{-1}BT^{-1} & = &
\left( \begin{array}{ccc}
\lambda_1\det t_1^{-1} J &  &  \raiseghost{-2ex}{-1ex}{{\LARGE 0}}\\
  &\ddots & \\
\raiseghost{2ex}{0ex}{{\LARGE 0}} & &\lambda_n\det t_n^{-1} J
\end{array} \right)
\end{array}
\]
this is equivalent to the question if there exist 
$\alpha,\beta,\gamma,\delta\in\CC$ with
\[
\begin{array}{c@{\;}c@{\;}c}
(\det t_1^{-1},\ldots,\det t_n^{-1}) & = & 
\alpha(1,\ldots,1)+\beta(\lambda_1,\ldots,\lambda_n)\\[1ex]
(\lambda_1 \det t_1^{-1},\ldots,\lambda_n \det t_n^{-1}) & = & 
\gamma(1,\ldots,1)+\delta(\lambda_1,\ldots,\lambda_n).
\end{array}
\]
It follows
\[
\begin{array}{l}
-\gamma(1,\ldots,1)+(\alpha-\delta)(\lambda_1,\ldots,\lambda_n)+
\beta(\lambda_1^2,\ldots,\lambda_n^2)=0\\[1ex]
\Longrightarrow \alpha =\delta,\ \beta=\gamma=0\\[1ex]
\Longrightarrow \det t_1=\ldots=\det t_n.  
\end{array}
\]
We normalize by $\det t_1=1$, i.e.\ $t_1,\ldots,t_n\in \SL(2,\CC)$.
Then only $T$ and $-T\in\GL(2n,\CC)$ give the same element 
in $\PP \GL(2n,\CC)$. 
So that as 
a group the automorphisms of $\GG(1,2n-1)\cap L$ that do not permute the 
exceptional lines are isomorphic to $\SL(2,\CC)^n/\{1,-1\}$.

\medskip
Altogether we get

\begin{theorem}\label{autog12n-1h2}
For $n\ge 3$ the automorphism group of the intersection of $\GG(1,2n-1)$ with
a general 2-codimensional linear subspace of $\PP(\bigwedge^2\CC^{2n})$ has
$\SL(2,\CC)^n/\{1,-1\}$ as a normal subgroup and 
the quotient group is isomorphic 
to the permutation group $\mathrm{S}(3)$ for $n=3$, 
to $\ZZ/2\ZZ\times\ZZ/2\ZZ$ for $n=4$, and trivial otherwise.

\medskip
The automorphism group is isomorphic to the subgroup of $\PP\GL(2n,\CC)$ 
that consists of the elements
\[
P_\sigma \cdot
\left( \begin{array}{ccc}
t_1 &  & \raiseghost{-2ex}{-1ex}{\emph{\LARGE 0}}\\
  &\ddots & \\
\raiseghost{2ex}{0ex}{\emph{\LARGE 0}} & &t_n
\end{array} \right)
\ \ \mbox{with\ }t_1,\ldots,t_n\in\SL(2,\CC)
\]
where $P_\sigma$ is the identity for $n\ge 5$ and otherwise defined by
\[
\begin{array}{c@{\;}c@{\;}l}
P_\sigma(e_{2i}) & = & e_{2\sigma(i)}\\
P_\sigma(e_{2i-1}) & = & e_{2\sigma(i)-1}
\end{array}
\]
\[ \mbox{for} \  \sigma \in \left\{
\begin{array}{ll}
\mathrm{S}(n) & \mbox{if\ }n=3\\[0.5ex]
\{ (1\,2\,3\,4),(2\,1\,4\,3),(3\,4\,1\,2),(4\,3\,2\,1)\}
& \mbox{if\ }n=4.
\end{array}\right.
\]
\end{theorem}

\bigskip
For the sake of completeness we recall the classical case of $\GG(1,3)\cap H^2$.

\begin{remark} The automorphism group of $\GG(1,3)\cap H^2$ is an extension of
$\ZZ/2\ZZ$ by $\PP\GL(2,\CC)\times\PP\GL(2,\CC)$. It acts homogeneously on 
$\GG(1,3)\cap H^2$.
\end{remark}
\emph{Proof.} The Grassmannian $\GG(1,3)$ is a smooth quadric in 
$\PP(\bigwedge^2\CC^4)\cong\PP_5$. Therefore $\GG(1,3)\cap H^2$ is a smooth 
quadric in $\PP_3$. Hence it is isomorphic to the Segre variety 
$\PP_1\times\PP_1$ in $\PP_3$. The automorphism group of $\PP_1\times\PP_1$ is
generated by $\PP\GL(2,\CC)\times\PP\GL(2,\CC)$ together 
with the automorphism
that
exchanges the $\PP_1$s. All the automorphisms extend to $\PP_3$. 
Obviously, the group acts transitively on $\PP_1\times\PP_1$.
\hfill$\Box$

\bigskip
For the rest of this section we consider the question if the action of the other 
automorphism groups is quasihomogeous on the corresponding line system, 
i.e.\ if there is an open orbit.

\medskip
This cannot be the case for $n\ge 7$ since then the dimension of the line system
$\GG(1,2n-1)\cap H^2$, $2(2n-2)-2=4n-6$, is larger than the dimension of the 
automorphism group, $3n$.

\medskip
For $n=3$ the action is quasihomogeneous. To see that one can adjust the 
$(\lambda_1,\lambda_2,\lambda_3)$ in the normal form of the line system to
$(1,0,-1)$ by a projective transformation and compute the stabiliser of the 
line $(1\,\colon 0\,\colon 1\,\colon 0\,\colon 1\,\colon 0)\wedge(1\,\colon 1\,\colon 1\,\colon -2\,\colon 1\,\colon 1)$ by hand or computer and see that it is
3-dimensional. So the dimension of its orbit is  $3\cdot 3-3=6$, which is 
just the dimension of the line system.

\medskip
For $n=4,5,6$ the group does not act quasihomogenously anymore. For this one 
computes again the dimension of the stabiliser of a  general line. Since the 
group acts transitively on $\PP_{2n-1}\setminus \bigcup L_i$, we may restrict
our attention to lines through one of those points, e.g. 
$(1\,\colon 0\,\colon\, \ldots\,\colon 1\,\colon 0)$.
Using a computer one sees that the stabilizer of a general line through this 
point has again dimension $3$. Hence the orbit has dimension $3n-3$,
 which is less then the dimension of the line system, $4n-6$.

\section{$\GG(1,5)\cap H^3$}

Let $L=H^3\subset \PP(\bigwedge^2\CC^6)$ be a general 3-codimensional subspace.
With our usual identification of $(\bigwedge^2\CC^6)^*$ with the antisymmetric
matrices $\antisym(6,\CC)$ its dual plane 
$L^*=\PP(\lambda A+\mu B +v C)\subset \PP(\bigwedge^2\CC^6)^*$
intersects the dual Grassmannian $\GG(1,5)^*$, which consists of matrices of 
rank $\le 4$ and is a hypersurface of degree 3 by Corollary \ref{dualG}, in an 
irreducible cubic $C^*$. By Corollary \ref{hyperplanetangency} a point
$(\lambda\,\colon \mu\,\colon \nu)\in C^*$  corresponds to the hyperplane
$h_{(\lambda\,\colon \mu\,\colon \nu)}=\PP(\lambda A+\mu B+\nu C)$ that is 
tangent to the Grassmannian at the point 
\[
l_{(\lambda:\mu:\nu)}:=\ker (\lambda A+\mu B+\nu C)\subset\PP_5.
\]
In analogy to the former case we have

\begin{lemma}
\[
\left\{ 
p\in\PP_5 \left|
\begin{array}{l}
\mbox{through}\ p\ \mbox{passes\ a\ }\PP_2 \\
\mbox{of}\ \mbox{lines\ of}\ \GG(1,5)\cap L
\end{array}
\right.\right\}
=\bigcup_{(\lambda:\mu:\nu)\in C^*} l_{(\lambda:\mu:\nu)} 
\subset\PP_5\]
\end{lemma}
\emph{Proof.} Since by definition
the lines in $\GG(1,5)\cap L$ that contain $p$ are $p\wedge q$ with 
${}^t\!pAq={}^t\!pBq={}^t\!pCq=0$, we see that 
\[
\begin{array}{l}
\mbox{through}\ p\ \mbox{passes\ at\ least\ a\ }\PP_2\
\mbox{of\ lines\ of\ } \GG(1,5)\cap L\\[0.5ex]
\Longleftrightarrow {}^t\!pA, {}^t\!pB,{}^t\!pC \ \mbox{are\ linear\ dependent}\\[0.5ex]
\Longleftrightarrow \exists {(\lambda\,\colon \mu\,\colon \nu)}\in\PP_2 \ \mbox{with}\ 
{}^t\!p(\lambda A+\mu B+\nu C)=0\\[0.5ex]
\Longleftrightarrow p\in\ker (\lambda A+\mu B+\nu C)=l_{(\lambda:\mu:\nu)}.
\end{array}
\]

We also note that there cannot be a $\PP_3$ of lines of $\GG(1,5)\cap L$ 
through a point $p$. Because if there were one, then 
$\dim \Span\{ {}^t\!p A,{}^t\!p B,{}^t\!pC\}=1$, i.e.\ there exist two points
$(\lambda\,\colon \mu\,\colon \nu),(\lambda'\,\colon \mu'\,\colon \nu')\in\PP_2$ with 
\[
{}^t\!p (\lambda A+\mu B+\nu C)={}^t\!p (\lambda' A+\mu' B+\nu' C)=0.
\]
It follows that all the matrices 
\[
(\alpha\lambda+\beta\lambda') A+(\alpha\mu+\beta\mu') B+(\alpha\nu+\beta\nu')C
\  \ \ \mbox{for\ all}\ (\alpha\,\colon \beta)\in\PP_1
\]
have a non-trivial kernel. Hence the line
$(\alpha\lambda+\beta\lambda'\,\colon \alpha\mu+\beta\mu'\,\colon \alpha\nu+\beta\nu')$
must lie in $L^*\cap \GG(1,5)^*=C^*$. But this is a contradiction 
since the cubic $C^*$ is irreducible.
\hfill$\Box$

\begin{proposition}
The lines $l_{(\lambda:\mu:\nu)}\subset\PP_5$ with 
$(\lambda\,\colon \mu\,\colon \nu)\in C^*$
do not intersect each other.
\end{proposition}
\emph{Proof.} Assume that the line $l_{(\lambda:\mu:\nu)}$ 
intersects the line
$l_{(\lambda':\mu':\nu')}$ in the point $p$, i.e.
\[
p\in\ker  (\lambda A+\mu B+\nu C) \cap \ker  
(\lambda' A+\mu' B+\nu' C)\not= 0.
\]
Then
\[
p\in \ker((\alpha\lambda+\beta\lambda') A+(\alpha\mu+\beta\mu') B
+(\alpha\nu+\beta\nu')C)\not= 0 \ \ \ 
\mbox{for\ all}\ (\alpha\,\colon \beta)\in\PP_1,
\]
and the line 
$(\alpha\lambda+\beta\lambda'\,\colon \alpha\mu+\beta\mu'\,\colon \alpha\nu+\beta\nu')$
must be  contained in the irreducible cubic $C^*$, which is a contradiction.
\hfill$\Box$

\bigskip
Let us again derive a normal form:

\begin{proposition}\label{normalcubic}
For a general plane 
$L^*=\PP (\lambda A+\mu B+\nu C)\subset\PP(\bigwedge^2\CC^6)^*$ 
there exists a choice of bases of $L^*$ and $\CC^6$ such that
\[
A=
\left(\begin{array}{cccccc}
0 &\!\! -1 &  & & &\raiseghost{-4ex}{-3.5ex}{\emph{\Huge 0}} \\
1 &  0 & & & & \\
 & &  0& \!\!-1 & & \\
 & &  1& 0 & & \\
 & &   &   &0 &0\\
\raiseghost{4ex}{1ex}{\emph{\Huge 0}} & & & &0 & 0
\end{array}\right)
\ \
B=
\left(\begin{array}{cccccc}
0 &  0 &  & & &\raiseghost{-4ex}{-3.5ex}{\emph{\Huge 0}} \\
0 &  0 & & & & \\
 & &  0& \!\!-1 & & \\
 & &  1& 0 & & \\
 & &   &   &0 &\!\!-1\\
\raiseghost{4ex}{1ex}{\emph{\Huge 0}} & & & &1 & 0
\end{array}\right)
\]
\[
C=
\left(\begin{array}{cc|cc|cc}
0 &0 & \!\!-\alpha  &0 & \!\! -\gamma&0\\
0 & 0 & 0&\!\!\!\!-\alpha &  \!\! -\delta&\!\!\!\!-\gamma\\ \hline
\alpha &0 & 0&\!\!\!\!-1 & \!\! -\beta &0\\
0&\alpha & 1&0 &0&\!\!\!\!-\beta\\ \hline
\gamma&\delta &\beta &0 & 0 &0 \\
0&\gamma &0&\beta & 0 & 0
\end{array}\right).
\]
\end{proposition}

\begin{remark} It is also possible to derive a more symmetric normal form
where all three matrices look like $C$ only with the ${ 0\,-1 \atop 1\ \;0}$
block moved along the diagonal, but this is not more useful for our 
computations.
\end{remark}
\emph{Proof of proposition \ref{normalcubic}.}
We may assume that the line $\PP(\lambda A+\mu B)\subset L^*$ is a
general line.
By Proposition \ref{normalpencil} there exists a choice of coordinates 
(corresponding to $\lambda_1=0$, $\lambda_2=1$, $\lambda_3=\infty$) such that
$A$ and $B$ are of the required form. Further, if we change the coordinates of 
$\CC ^6$ by transformations of the type
\[
T=
\left(\begin{array}{ccc}
t_1 & 0 & 0\\
0 & t_2 & 0 \\
0 & 0 & t_3
\end{array}\right)
\ \ t_1,t_2,t_3\in \SL(2,\CC),
\]
then $A$ and $B$ will stay the same by Theorem \ref{autog12n-1h2}.

\medskip
We write the matrix $C$ as
\[
C=\left(\begin{array}{ccc}
c_1 J &\!\! -^tC_{21} &\!\! -^tC_{31}\\[0.5ex]
C_{21} & c_2J &\!\! -^t C_{32} \\[0.5ex]
C_{31} & C_{32} & c_3 J
\end{array}\right)
\ \ 
\mbox{with}\ \ 
\begin{array}{l}
J=\left(\begin{array}{cc}0&-1\\1&0\end{array}\right);
\ c_1,c_2,c_3\in\CC\\[2.5ex]
C_{21},C_{31}, C_{32}\in \M(2,\CC).
\end{array}
\]
We may assume that $c_1=c_3=0$, $c_2=1$, otherwise we replace $C$ by the 
matrix
$1/(c_2-c_1-c_3)(C-c_1 A-c_3B)$. This is possible since $c_2-c_1-c_3\not=0$,
because $C$ is general. So C looks like
\[
C=\left(\begin{array}{ccc}
0 &\!\! -^tC_{21} &\!\! -^tC_{31}\\[0.5ex]
C_{21} & J &\!\! -^t C_{32} \\[0.5ex]
C_{31} & C_{32} & 0
\end{array}\right).
\]
The generality of $C$ ensures that the matrices $C_{21}$ and $C_{32}$
are invertible, so
\[
T=\left(\begin{array}{ccc}
\frac{1}{\alpha} C_{21} & 0 &0\\[0.5ex]
0 & \E_2 & 0 \\[0.5ex]
0 & 0 &\frac{1}{\beta} {}^t\!C_{32}
\end{array}\right)
\ \ 
\mbox{with}\ \ 
\begin{array}{l}
\alpha=\sqrt{\det C_{21}}\\[0.5ex]
\beta=\sqrt{\det C_{32}}
\end{array}
\]
is of the above mentioned type and transforms $C$ into
\[
C':={}^tT^{-1}CT^{-1}=\left(\begin{array}{ccc}
0 &\!\! -\alpha \E_2 &\!\! -^t\overline{C}\\
\alpha\E_2 & J &\!\! -\beta \E_2 \\
\overline{C} &\beta \E_2 & 0
\end{array}\right)
\ \ \mbox{with} \ 
\overline{C}:=\alpha\beta C_{32}^{-1}C_{31}C_{21}^{-1}.
\]
This matrix will be transformed under 
\[
T=\left(\begin{array}{ccc}
t^{-1}\!\!\! & 0 &0\\
0 & {}^t t & 0 \\
0 & 0 & t^{-1}\!\!\!
\end{array}\right)
\ \ 
\mbox{with}\ t\in\SL(2,\CC)
\]
into
\[
{}^tT^{-1}C'T^{-1}=\left(\begin{array}{ccc}
0 & \!\!-\alpha \E_2 & \!\!-{}^t t ^t\overline{C}t\\[0.5ex]
\alpha\E_2 & J & \!\!-\beta \E_2 \\[0.5ex]
{}^t t\overline{C}t &\beta \E_2 & 0
\end{array}\right).
\]
So, all that remains to show is: Given a general matrix 
$\overline{C}\in\M(2,\CC)$ there is a matrix $t\in\SL(2,\CC)$ such that
\[
{}^t t\overline{C}t=
\left(\begin{array}{cc}
\gamma & \delta\\
0 & \gamma 
\end{array}\right).
\]

If 
\[
\overline{C}=\left(\begin{array}{cc}
c_{11} & c_{12} \\[0.5ex]
c_{21} & c_{22}
\end{array}\right) 
\ \ \mbox{and} \ \ 
t=
\left(\begin{array}{cc}
1 &\!\! -\frac{c_{21}}{c_{11}}\\[0.5ex]
0 & 1
\end{array}\right) 
\]
then
\[
\overline{C}'={}^t t\overline{C}t=
\left(\begin{array}{cc}
c_{11} & c_{12} -c_{21}\\[0.5ex] 
0 & \frac{\det \overline{C}}{c_{11}}
\end{array}\right) 
\]
and an additional transformation by
\[
t=
\left(\begin{array}{cc}
 \frac{\sqrt[4]{\det \overline{C}}}{\sqrt{c_{11}}} & 0\\[0.5ex]
0 &  \frac{\sqrt{c_{11}}}{\sqrt[4]{\det \overline{C}}}
\end{array}\right) 
\]
takes $\overline{C}$ into the desired form

\bigskip
\ \hfill 
$\displaystyle 
{}^t t\overline{C}'t=
\left(\begin{array}{cc}
\sqrt{\det \overline{C}} & c_{12} -c_{21} \\[0.5ex]
0 & \sqrt{\det \overline{C}}
\end{array}\right) .
$
\hfill$\Box$

\bigskip
\begin{remark}
In terms of this coordinates the cubic $C^*\subset L^*$ is given as
\[
\lambda^2\mu+\mu^2\lambda+\lambda\mu\nu-(\gamma^2+\beta^2)\lambda\nu^2
-(\alpha^2+\gamma^2)\mu\nu^2+(\alpha\beta\delta -\gamma^2)\nu^3.
\]
One checks that the cubic is smooth for general $\alpha,\beta,\gamma,\delta$.
\end{remark}

\bigskip
Now we start to determine the automorphism group of $\GG(1,5)\cap L$.
A given automorphism
\[
\phi\in\aut(\GG(1,5)\cap L)=
\aut(L,\PP({\textstyle \bigwedge\nolimits^2}\CC^6))\cap
\aut(\GG(1,5),\PP({\textstyle \bigwedge\nolimits^2}\CC^6))
\]
induces a dual automorphism $\phi^*$ on the dual projective space
$\PP(\bigwedge\nolimits^2\CC^6))^*$
that preserves $L^*$ and the dual 
Grassmannian $\GG(1,5)^*$, i.e.
\[
 \phi^*\in
\aut(L^*,\PP({\textstyle \bigwedge\nolimits^2}\CC^6)^*)\cap
\aut(\GG(1,5)^*,\PP({\textstyle \bigwedge\nolimits^2}\CC^6)^*).
\]
In particular, $\phi^*$ induces a projective transformation of $L^*$ 
preserving $C^*$. But a smooth cubic has only finitely many
 automorphisms that are 
induced by a projective linear transformation \cite[7.3]{BK}. 

\bigskip
To find  all automorphisms of 
$\GG(1,5)\cap L$ that induce the identity on $L^*$, 
we look for the 
$T\in\PP\GL(6,\CC)$ such that
\[
{}^tT^{-1}(\lambda A+\mu B+\nu C)T^{-1} \in \CC \cdot (\lambda A+\mu B+\nu C)
\ \ \ \mbox{for\ all}\ \lambda,\mu,\nu\in\CC.
\]

% we look in fact for $T$ such that
%\[
%{}^tT^{-1}(\lambda A+\mu B+\nu C)T^{-1} = \lambda A+\mu B+\nu C
%\ \ \ \mbox{for\ all}\ \lambda,\mu,\nu\in\CC.
%\]
It suffices to check this for $(\lambda,\mu,\nu)=(1,0,0)$, $(0,1,0)$,
and $(0,0,1)$.
If we normalize the representation of $T$ in $\GL(6,\CC)$ by $\det T=1$,
we know from the previous section that ${}^tT^{-1}AT^{-1}=\CC \cdot A$ 
and ${}^tT^{-1}BT=\CC \cdot B$ is equivalent to 
\[
T=
\left(\begin{array}{ccc}
t_1 & 0 & 0\\
0 & t_2 & 0 \\
0 & 0 & t_3 
\end{array}\right)
\ \ \mbox{with}\ t_1,t_2,t_3\in\SL(2,\CC).
\]
Furthermore, we compute 
\[
{}^tT^{-1}CT^{-1}=
\left(\begin{array}{ccc}
0 & -\alpha\;\! {}^t t_1^{-1}t_2^{-1} &
\!\! -{}^t t_1^{-1}{\gamma\, 0 \choose \delta\, \gamma}t_3^{-1} \\[2ex]
\alpha\;\! {}^t t_2^{-1}t_1^{-1} & 0 &  -\beta\;\! {}^t t_2^{-1}t_3^{-1} \\[2ex]
{}^t t_3^{-1}{\gamma\, \delta \choose 0 \, \gamma}t_1^{-1} &  \beta\;\! {}^t t_3^{-1}t_2^{-1} & 0
\end{array}\right),
\]
so that ${}^tT^{-1}CT^{-1}=\vartheta \cdot C$ iff $t_1=\frac{1}{\vartheta} 
{}^t t_2^{-1}=t_3=:t$ and
\[
{}^t t^{-1}\left(\begin{array}{cc}
\gamma & \delta\\
0 & \gamma
\end{array}\right)t^{-1}
= \vartheta
\left(\begin{array}{cc}
\gamma & \delta\\
0 & \gamma
\end{array}\right).
\]
Because of $\det t_1=\det t_2=1$, $\vartheta$ must be either 1 or -1. Setting
\[
t=
\left(\begin{array}{cc}
a & b\\
c & d
\end{array}\right)
\Longrightarrow 
t^{-1}=
\left(\begin{array}{cc}
d & -b\\
-c & a
\end{array}\right)
\]
the last condition together with $\det t=1$ requires that the following 
polynomials vanish:
\[
\begin{array}{l}
(d^2+c^2-\vartheta)\gamma- dc \delta,\ (db+ac)\gamma-(\vartheta-ad)\delta\\[0.5ex]
(db+ac)\gamma-bc\delta,\ (b^2+a^2-\vartheta)\gamma-ba\delta,\ ad-bc-1 
\end{array}
\]

The Gr\"obner basis of the ideal generated by these polynomials with
respect to the lexicographical order $\gamma>\delta>a>b>c>d$ can be computed 
for $\vartheta=1$ as
\[
\gamma a+\delta c -\gamma d,\  b+c ,\  ad+c^2-1,
\]
so that 
\[
t=
\left(\begin{array}{cc}
a&\!\frac{\gamma}{\delta}(a-d)\\[0.75ex]
\!\frac{\gamma}{\delta}(d-a) & d
\end{array}\right)
\ \ \mbox{with} \ \det t=1.
\]
For $\vartheta=-1$ we get as the Gr\"obner basis 
\[
\delta, a+d,-c+b,d^2+c^2+1.
\]
Since in the general case $\delta\not=0$, this gives no further automorphisms.

\medskip
The one-dimensional subgroup of $\PP\GL(2,\CC)$ consisting of elements like 
$t$ above acts on $\PP_1$ with the two fixed points 
$(-\delta \pm \sqrt{\delta^2-4\gamma^2}:2 \gamma)$. Hence it is conjugate to
the one-dimensional subgroup of $\PP\GL(2,\CC)$ that acts on $\PP_1$ with
 the fixed points 0 and $\infty$. Now this subgroup consists of the 
invertible diagonal matrices of $\PP\GL(2,\CC)$, so it is isomorphic to
$\CC^*$. Therefore we have shown 

\begin{theorem}
The component of the automorphism group of $\GG(1,5)\cap H^3$
containing the identity is isomorphic to $\CC^*$.
The quotient of $\aut(\GG(1,5)\cap  H^3)$ by this component is a subgroup
of the finite group of projective automorphisms of a smooth cubic in 
$\PP_2$.
\end{theorem}

\section{$\GG(1,2n)\cap H$}

The hyperplane $H$ is given by an element $A\in (\bigwedge^2\CC^{2n+1})^*$
which can be thought of as an antisymmetric matrix of size $2n+1$. Since
antisymmetric matrices have an even rank, the general $H$ corresponds to 
an $A$ of rank $2n$. The one dimensional kernel of $A$ as a point of $\PP_{2n}$
is called the center $c$ of $H$.

\medskip
The center plays a special role in the geometry of the line system 
$\GG(1,2n)\cap H$ in $\PP_{2n}$.

\begin{proposition}\label{centerprop}
Every line through the center of the line system $\GG(1,2n)\cap H$ is in the 
line system. The center is the only point with this property.

\medskip
Moreover, if the line $l\not\ni c$ belongs to the line system, so does every 
line in the plane spanned by the line $l$ and the center $c$.
\end{proposition}
\emph{Proof.} The line $c\wedge p$ through the center will be in the 
line system if ${}^t\!c A p=0$. But $c$ is the kernel of $A$, so this is true.
On the other hand, if $\overline{c}$ is a point such that every line 
through it belongs to the line system, then ${}^t\overline{c} A p=0$ for
all $p\in \PP_{2n}$. Hence  $\overline{c}$ must be in the kernel of $A$, 
and therefore  $\overline{c}=c$.

\medskip
Let the line $l=p\wedge q$ be in the line system. All the lines in the plane
spanned by $l$ and $c$ -- except the lines through $c$ itself --  
can be written as
\[
(\alpha p +\beta c)\wedge (\lambda q+\mu c) \quad
\mbox{for\ }(\alpha:\beta),(\lambda:\mu)\in \PP_1.
\]
These will be in the line system since

\bigskip
\ \hfill $\displaystyle
\begin{array}{c@{\;=\;}l}
(\alpha {}^t\!p +\beta {}^t\!c) A (\lambda q+\mu c) &
\alpha\lambda {}^t\!p A q+\alpha\mu {}^t\!p A c
+\beta\lambda  {}^t\!c A q+\beta\mu  {}^t\!c A c\\[0.5ex]
& \alpha\lambda {}^t\!p A q =0.
\end{array}
$ 
\hfill$\Box$

\bigskip
Let us for a moment look at the projection $\PP(\CC^{2n+1}/c)$ of 
$\PP_{2n}$ from the center $c$. This projection maps all lines in a plane 
through the center -- except the lines through $c$ itself -- to only one line.
Hence we get a codimension one line system inside $\PP_{2n-1}$. In fact, it is 
of the form $\GG(1,2n-1)\cap \overline{H}$, which is most easily seen in 
coordinates. We choose a basis $(e_0,\ldots,e_{2n})$ of $\CC^{2n+1}$ such that
the hyperplane $H$ is given by the matrix
\[
A=\left(
\begin{array}{c|c}
\begin{array}{cc}
{\displaystyle 0} &{\displaystyle -\E_n}\\[0.5ex]
{\displaystyle\E_n} &{\displaystyle 0}
\end{array} 
&
\begin{array}{c} \!\! 0 \\[-0.8ex]\!\! \vdots  \\[-0.3ex] \!\!0 \end{array}
\\ \hline
0\, \cdots\, 0 
&\!\! 0
\end{array}
\right)\in \antisym(2n+1,\CC).
\]
The center of $H$ is $c=\PP(e_{2n})$. So the projected line system is 
$\GG(1,2n-1)\cap \overline{H}$, where $\overline{H}$ is given by the matrix $A$
with the last row and column deleted.

\bigskip
This description helps to determine the automorphism group of 
$\GG(1,2n) \cap H$.

\medskip
First of all, any of the automorphism must -- as a transformation 
$T\in\PP\GL (2n+1,\CC)$ -- preserve the center, i.e. $Tc=c$. Therefore it
induces a transformation $\overline{T}$ of the projected space 
$\PP(\CC^{2n+1}/c)$. This induced transformation  $\overline{T}$ has to 
preserve the projected line system $\GG(1,2n-1)\cap \overline{H}$.
Since this case has been treated in Section
\ref{g12n-1hsection}, we know that if we  normalize  $\overline{T}$ 
by $\det \overline{T}=1$, then  $\overline{T}\in \Sp(2n,\CC)$. Therefore
$T$ must have been of the form 
\[
T=
\left(
\begin{array}{c|c}
{\displaystyle \overline{T}}
&
\begin{array}{c}
\!\!0 \\[-0.8ex]\!\! \vdots  \\[-0.3ex]\!\! 0
\end{array}
\\ \hline
a_0 \, \cdots\, a_{2n-1} 
& \!\! b
\end{array}
\right)
\quad \mbox{with}\ 
\begin{array}{l}
\overline{T}\in\Sp(2n,\CC)\\
a_i\in \CC \\
b\in \CC^*.
\end{array}
\]

One immediately checks that ${}^tT^{-1}AT=A$, 
so that the automorphism group as a 
subset of $\PP\GL(2n+1,\CC)$ consists of all elements of the above type.
Since we normalized $\overline{T}$, we have up to multiplication by $-1$ 
an unique representative in the class of $\PP\GL(2n+1,\CC)$.

\medskip
A small computation shows that 
\[
N:=\left\{
\left(
\begin{array}{c|c}
{\displaystyle \E_{2n}}
&
\begin{array}{c}
\!\!0 \\[-0.8ex]\!\! \vdots  \\[-0.3ex]\!\! 0
\end{array}
\\ \hline
a_0 \, \cdots\, a_{2n-1} 
& \!\! 1
\end{array}
\right)
\mid a_i\in\CC
\right\} \subset \aut(\GG(1,2n)\cap H)
\]
is a normal subgroup  which is isomorphic to $(\CC^{2n},+)$.

\medskip
Collecting everything together we have

\begin{proposition} The automorphism group of $\GG(1,2n)\cap H$ for a general
hyperplane $H\subset \PP(\bigwedge^2\CC^{2n+1})$ is an extention of
$\Sp(2n,\CC)\times\CC^*/\{1,-1\}$ by $(\CC^{2n},+)$ and is isomorphic 
to the group
\[
\left\{
\left(
\begin{array}{c|c}
{\displaystyle T}
&
\begin{array}{c}
\!\!0 \\[-0.8ex]\!\! \vdots  \\[-0.3ex]\!\! 0
\end{array}
\\ \hline
a_0 \, \cdots\, a_{2n-1} 
& \!\! b
\end{array}
\right) 
\begin{array}{|l}
T\in\Sp(2n,\CC)\\
a_i\in \CC \\
b\in \CC^*
\end{array}
\right\}/\{1,-1\}.
\]
\end{proposition}

\bigskip
The action of the automorphism group on the line system is described by the
following

\begin{proposition} The action of the automorphism group of $\GG(1,2n)\cap H$
on the lines of $\GG(1,2n)\cap H$ has two orbits:
\begin{enumerate}
\item the lines containing the center $c$
\item the lines that do not.
\end{enumerate}
\end{proposition}
\emph{Proof.} Since all the automorphisms preserve the center any orbit 
will be contained in these two sets. 

\medskip
First we show that the lines containing
$c$ form one orbit. For two lines $c\wedge p$ and $c\wedge q$, we may 
assume $p,q\in\PP(\CC^{2n}\!\times 0)$. Take a $\overline{T}\in \Sp(2n,\CC)$
that maps $p$ to $q$. The trivial extention of $\overline{T}$ to 
$T\in\SL(2n+1,\CC)$ will take $c\wedge p$ to $c\wedge q$.

\medskip
The other lines will form the second orbit since any line not containing
the center can be pushed into the hyperplane $\PP(\CC^{2n}\!\times 0)$ by 
a transformation with an element of the normal subgroup $N$. There one can
use the transitive action of the $\aut(\GG(1,2n-1)\cap H)$ subgroup to
show that all these lines can be mapped onto each other.
\hfill$\Box$

\section{$\GG(1,2n)\cap H^2$}

Let $L=H^2$ be a 2-codimensional linear subspace of 
$\PP(\bigwedge^2\CC^{2n+1})$. We want to study the linear line system
$\GG(1,2n)\cap L$. To $L$ corresponds the line 
$L^*=\PP(\lambda A-\mu B)\subset \PP( \bigwedge^2\CC^{2n+1} )^*$ of the 
hyperplanes $H_{(\lambda:\mu)}=\PP(\lambda A-\mu B)$ containing $L$. We 
identify as always $(\bigwedge^2\CC^{2n+1})^*$ with the antisymmetric matrices
$\antisym(2n+1,\CC)$. The locus of antisymmetric matrices of corank~3 in
$\antisym(2n+1,\CC)$ is 3-codimensional by Corollary \ref{dualG}. 
Therefore a line $L^*$ may 
be called general if it does not intersect it. Hence for the general 
line $L^*$ the antisymmetric matrices $\lambda A-\mu B$ corresponding to 
the hyperplanes $H_{(\lambda:\mu)}$  have all corank~1. So each of the 
hyperplane sections  $\GG(1,2n)\cap H_{(\lambda:\mu)}$ has a unique 
center $c_{(\lambda:\mu)}\in \PP_{2n}$ by Proposition \ref{centerprop}.
These centers play an important role in the geometry of $\GG(1,2n)\cap L$.

\begin{proposition}
The centers $c_{(\lambda:\mu)}$ are those points of $\PP_{2n}$ through 
which there passes a $\PP_{2n-2}$ of lines of the line system 
$\GG(1,2n)\cap L$. Through all the other points of $\PP_{2n}$ passes only a
$\PP_{2n-3}$ of lines.
\end{proposition}
\emph{Proof.} The lines of the line system through a point $p\in \PP_{2n}$
are $p\wedge q$ with ${}^t\!p Aq={}^t\!pBq=0$. 
So we need to show that ${}^t\!p A$ and 
${}^t\!pB$ are linear dependent iff $p$ is a center of a hyperplane 
$H_{(\lambda:\mu)}$.
Now  ${}^t\!p A$ and ${}^t\!pB$ are linear dependent precisely if there exists a 
$(\lambda\,\colon \mu)\in\PP_1$ with 
$0=\lambda {}^t\!pA-\mu{}^t\!pB={}^t\!p (\lambda A-\mu B)$,
i.e.\ $p$ is the kernel of $\lambda A-\mu B$, which is by definition the center
of $H_{(\lambda:\mu)}$.
\hfill $\Box$

\bigskip
\begin{remark}\label{twocenteroncurve}
Any line that contains two centers is a member of the line system $\GG(1,N)\cap L$.
\end{remark}
\emph{Proof.} If the line contains the centers $c_{(\alpha:\beta)}$ and 
$c_{(\lambda:\mu)}$, it is contained in the hyperplanes $H_{(\alpha:\beta)}$ and 
$H_{(\lambda:\mu)}$ by Proposition \ref{centerprop} and therefore in their
intersection $L=H_{(\alpha:\beta)}\cap H_{(\lambda:\mu)}$.
\hfill$\Box$

\bigskip

Next we want to know more about the curve $c_{(\lambda:\mu)}$.

\begin{proposition} 
Let $A$, $B$ be two antisymmetric matrices of size $2n+1$ such that every 
non-zero linear combination of them has corank 1. Then the map
\[
\begin{array}{c@{\,}ccc}
c: & \PP_1 & \longrightarrow & \PP_{2n}\\[0.3ex]
 & (\lambda\,\colon \mu) &\longmapsto & \ker (\lambda A-\mu B)
\end{array} 
\]
is a parametrisation of a rational normal curve of degree $n$.
\end{proposition}
\emph{Proof.} (compare \cite[X,4.3]{SR} for $n=2$.) 
First we show that the map is 
injective. If it is not, there are two points of $\PP_1$ with the same image.
We may assume that this is the case for $(1\,\colon 0)$ and $(0\,\colon 1)$, i.e.\ $A$ and
$B$ have the same kernel, say $e_0$. Writing $A$ and $B$ in a basis with $e_0$
as first element, we have
\[
A=
\left(\!\!\!\begin{array}{cc}
0 &\!\!\!\!\!\!\!\; \cdots\, 0 \\[-0.8ex]
\begin{array}{c} \vdots \\[-0.2ex]
0 \end{array} & \!\!\!\!\!\!\!\; {\displaystyle \widetilde{A}}
\end{array} \right)
\ \mbox{and}\ \
B=
\left(\!\!\!\begin{array}{cc}
0 &\!\!\!\!\!\!\!\; \cdots\, 0 \\[-0.8ex]
\begin{array}{c} \vdots \\[-0.2ex]
0 \end{array} & \!\!\!\!\!\!\!\; {\displaystyle \widetilde{B}}
\end{array} \right)
\ \ \mbox{with}\ \widetilde{A},\widetilde{B}\in\antisym(2n,\CC).
\]

Since $\det(\lambda \widetilde{A}-\mu \widetilde{B})$ is a homogeneous 
polynomial of degree $2n$, there exist a $(\lambda'\,\colon \mu')\in\PP_1$ with
$\det(\lambda' \widetilde{A}-\mu' \widetilde{B})=0$. But then 
$\lambda' A-\mu' B$ has corank at least two, which contradicts our assumption.

\medskip
Secondly, we proof that the map is of maximal rank everywhere. If it is not, we
may assume that it is not maximal at $(1\,\colon 0)$. Restricting to the chart 
$\lambda=1$, this means $c'(0)=0$. Now from
\[
\begin{array}{l}
(A-\mu B)c(\mu)=0\\[0.8ex]
\Longrightarrow  Ac'(\mu) - B c(\mu)-\mu B c'(\mu)=0 \\[0.8ex]
\Longrightarrow  Ac'(0) - B c(0)=0 
\end{array}
\]
$Bc(0)=0$ follows. Therefore $A$ and $B$ have the same kernel $c(0)$, and we 
are back in the above chain of arguments.

\medskip
Finally, we have to show that the embedding $c$ is of degree $n$. For this we
give an explicit form of the map. Recall \cite[5.2]{B} that the determinant of an
antisymmetric matrix $C=(c_{ij})$ of size $2n$ is the square of the irreducible Pfaffian
polynomial $\mathrm{Pf}\, C$,
\[
\mathrm{Pf}\, C:=\sum_\sigma \sgn(\sigma) c_{\sigma(1) \sigma(2)}\ldots
 c_{\sigma(2n-1) \sigma(2n)},
\]
where $\sigma$ runs through all permutations $\mathrm{S}(2n)$ with
$\sigma(2i-1)<\sigma(2i)$ for $i=1\ldots n$ and $\sigma(2i)<\sigma(2i+2)$
for $i=1\ldots n-1$.

\medskip
Let $c_i(\lambda\,\colon \mu)$ denote $(-1)^i$-times the Pfaffian of the matrix
$\lambda A-\mu B$ with the $i$-th row and column deleted. Then the 
$c_i(\lambda\,\colon \mu)$ are irreducible polynomials of degree $n$ and by a
straightforward but messy computation one can check that 
$(c_0(\lambda\,\colon \mu)\,\colon \ldots\,\colon c_{2n}(\lambda\,\colon \mu))$
is the kernel of $\lambda A-\mu B$. Therefore $c=(c_0\,\colon \ldots\,\colon c_{2n})$, 
which
shows that $c$ is a degree $n$ embedding of $\PP_1$.
\hfill$\Box$

\bigskip
After we have determined the special points in $\PP_{2n}$ of the line system
$\GG(1,2n)\cap L$, we are nearly ready to compute its automorphism group.
It remains to give a normal form for the line 
$L^*\subset \PP(\bigwedge^2\CC^{2n+1})$ to make computations easier.
This normal form was found by Donagi \cite[2.2]{D}. But he did not give
a proof for it since his main interest was lines in $\PP_{2n-1}$ and
not in $\PP_{2n}$. So we give the proof here.

\begin{proposition}
Let $L^*$ be a line in $\PP(\bigwedge^2\CC^{2n+1})^*$ such that the 
antisymmetric
matrices corresponding to the points of $L^*$ have all corank 1. Then there 
exists a basis $(e_0,\ldots,e_{2n})$ of $\CC^{2n+1}$ such that the line can 
be taken as $L^*=\PP(\lambda A-\mu B)$ with the matrices
\[ 
A=
\left(\begin{array}{c|c}
\begin{array}{cc}
{\displaystyle 0}  &\!\! {\displaystyle-\E_n}\\[0.5ex]
{\displaystyle\E_n} &{\displaystyle 0}
\end{array} 
& 
\begin{array}{c}
\!\! 0 \\[-0.8ex]\!\! \vdots \\[-0.3ex]\!\! 0
\end{array} \\ \hline
0\, \cdots\,  0 & \!\! 0
\end{array}\right)
\ \ \mbox{and} \ \ \ 
B=
\left( \begin{array}{c|c|c}
{\displaystyle 0}\! & 
\begin{array}{c}\!\!\! 0\!\!\!\!  \\[-0.8ex] \!\!\!\vdots\!\!\!\!  \end{array} &
\! {\displaystyle -\E_n} \\ \cline{1-1}\cline{3-3}
\multicolumn{3}{c}{\hspace{0.4ex} 0\; \cdots \; 0 \;\cdots\; 0} \\[-0.35ex]\cline{1-1}\cline{3-3}
{\displaystyle \raiseghost{4ex}{0ex}{${\displaystyle \E_n}$}\phantom{-\E_n}}\! & 
\begin{array}{c} \!\!\!\vdots\!\!\!\!  \\[-0.3ex]\!\!\!0\!\!\!\! \end{array} 
&\! {\displaystyle 0} 
\end{array}\right).
\]
\end{proposition}
\emph{Proof.} Let $A$ and $B$ any two matrices of $L^*$ in an arbitrary basis.
We will adjust the basis in three steps to achieve the required form for 
$A$ and $B$.

\medskip
\noindent\emph{$1^{st}$ Step:}
We know that the map
\[
\begin{array}{c@{\,}ccc}
c: & \PP_1 & \longrightarrow & \PP_{2n}\\[0.5ex]
 & (\lambda\,\colon \mu) &\longmapsto & \ker (\lambda A-\mu B)
\end{array} 
\]
is a parametrisation of a rational normal curve of degree $n$. Modulo 
projective transformations of $\PP_1$ and $\PP_{2n}$ such parametrisations are
all the same. So we can pick a basis of $\PP_1$ and $n+1$ linear independent
vectors $e_n,\ldots,e_{2n}$ of $\CC^{2n+1}$ such that 
\[
\begin{array}{c@{\,}ccc}
c: & \PP_1 & \longrightarrow & \PP_{2n}\\[0.5ex]
 & (\lambda\,\colon \mu) &\longmapsto & 
\PP \left(\sum\limits_{i=0}^n \lambda^i\mu^{n-i} e_{n+i}    \right).
\end{array} 
\]

Extending $(e_n,\ldots,e_{2n})$ to a basis $(e_0,\ldots,e_{2n})$ of
$\CC^{2n+1}$ and denoting by $a_0,\ldots,a_{2n}$ resp. $b_0,\ldots,b_{2n}$
the columns of $A$ resp. $B$,
the fact that $c(\lambda\,\colon \mu)$ is the kernel of
$\lambda A-\mu B$ for all $(\lambda\,\colon \mu)\in\PP_1$ has the following consequences 
for $A$ and $B$:

\bigskip
$ \hspace{9ex}( \lambda A-\mu B) 
\left(\sum\limits_{i=0}^n \lambda^i\mu^{n-i} e_{n+i} \right)=0$\vspace{1ex}  

$\hspace{9ex} \Longrightarrow \sum\limits_{i=0}^n \lambda^{i+1}\mu^{n-i} a_{n+i}
-\sum\limits_{i=0}^n \lambda^i\mu^{n+1-i} b_{n+i} =0  $\vspace{1ex} 

$\hspace{9ex}\Longrightarrow -\mu^{n+1}b_n
+\sum\limits_{i=0}^{n-1} \lambda^{i+1}\mu^{n-i}(a_{n+1}-b_{n+i+1})
+\lambda^{n+1} a_{2n} =0 $\vspace{1ex} 

$\hspace{9ex}\Longrightarrow b_n=0,\ a_{2n}=0,\ a_{n+i}=b_{n+i+1}
\ \ \mbox{for\ }i= 0\ldots n-1. \hfill (*)$

\bigskip
We claim that this implies:
\[
a_{n+i,n+j}=0 \ \ \ \ \mbox{for\ }i,j=0\ldots n.
\]

Indeed, for $1\le i \le n$,  $0\le j\le n-1$ we have using $(*)$
\[
a_{n+i,n+j}=b_{n+i,n+j+1}=-b_{n+j+1,n+i}=-a_{n+j+1,n+i-1}=a_{n+i-1,n+j+1}.
\]
This shows that the $a_{n+i,n+j}$ are all the same for $i+j=const$, in 
particular $a_{n+i,n+j}=a_{n+j,n+i}$. On the other hand, by the 
antisymmetricity of $A$ we have  $a_{n+i,n+j}=-a_{n+j,n+i}$, and the claim 
follows.

\medskip
Using $(*)$ again we know that $A$ and $B$ look in our basis like
\[ 
A=
\left(\begin{array}{c|c}
\begin{array}{cc}
{\displaystyle \widetilde{A}}  &\!\! {\displaystyle-{}^t\! M}\\[0.5ex]
{\displaystyle M} &{\displaystyle 0}
\end{array} 
& 
\begin{array}{c}
\!\! 0 \\[-0.8ex]\!\! \vdots \\[-0.3ex]\!\! 0
\end{array} \\ \hline
0\, \cdots\,  0 & \!\! 0
\end{array}\right)
\ \ \mbox{and} \ \ 
B=
\left( \begin{array}{c|c|c}
{\displaystyle \widetilde{B}}\! & 
\begin{array}{c}\!\!\! 0\!\!\!\!  \\[-0.8ex] \!\!\!\vdots\!\!\!\!  \end{array} &
\! {\displaystyle -{}^t\! M} \\ \cline{1-1}\cline{3-3}
\multicolumn{3}{c}{\hspace{0.4ex} 0\; \cdots \; 0 \;\cdots\; 0} \\[-0.35ex]\cline{1-1}\cline{3-3}
{\displaystyle \raiseghost{4ex}{0ex}{${\displaystyle M}$}\phantom{-{}^t\! M}}\! & 
\begin{array}{c} \!\!\!\vdots\!\!\!\!  \\[-0.3ex]\!\!\!0\!\!\!\! \end{array} 
&\! {\displaystyle 0} 
\end{array}\right)
\]
with $\widetilde{A},\widetilde{B}\in \antisym(n,\CC)$ and $M\in\GL(n,\CC)$.

\medskip
\noindent\emph{$2^{nd}$ Step:}
Here we will improve the choice of $(e_0,\ldots,e_{n-1})$ to achieve 
$\widetilde{A}=0$ and $M=\E_n$. We claim:

\medskip
Let an antisymmetric matrix $A\in\antisym(2n+1,\CC)$ of rank $2n$ and linear
independent vectors $e_n,\ldots, e_{2n}\in\CC^{2n+1}$ with ${}^t\!e_iAe_j=0$ 
for
$n\le i,j\le 2n$ and $Ae_{2n}=0$ be given. 
Then $(e_n,\ldots,e_{2n})$ can be extended 
to a basis $(e_0,\ldots,e_{2n})$ of $\CC^{2n+1}$ 
such that in this basis $A$ is 
given as 
\[
A=\left(\begin{array}{c|c}
\begin{array}{cc}
{\displaystyle 0}  & \!\!{\displaystyle-\E_n}\\[0.5ex]
{\displaystyle\E_n} &{\displaystyle 0}
\end{array} 
& 
\begin{array}{c}
\!\! 0 \\[-0.8ex]\!\! \vdots \\[-0.3ex]\!\! 0
\end{array} \\ \hline
0\, \cdots\,  0 & \!\! 0
\end{array}\right) 
.\]

The proof is by induction. The statement is trivial for $n=0$. Assuming the claim 
for $n-1$, we prove it for $n$. 
Let $W:=\bigcap_{i=1}^n\ker {}^t\!e_{n+i}A$, then 
there exists an $e_0\in W$ with ${}^t\!e_n A e_0=1$. If not, we would have 
$W=W\cap \ker {}^t\!e_nA$ and with $e_{2n}\in \ker A$
\[
\dim\, \Span\{{}^t\!e_nA,\ldots,{}^t\!e_{2n}A\}=
\dim\, \Span\{{}^t\!e_{n+1} A,\ldots,{}^t\!e_{2n-1}A\}\le n-1,
\]
which contradicts $\rank A=2n$.

\medskip
Set $V:=\ker {}^t\!e_0A \cap \ker{}^t\!e_nA$, then $\dim V =2(n-1)+1$ and 
$e_{n+1},\ldots,e_{2n}\in V$. Therefore the induction hypothesis can be applied
to $\left.A\right|_V$. Together with ${}^t\!e_0Ae_n=1$ and 
${}^t\!e_0Av={}^t\!e_nAv=0$ for $v\in V$ this implies the stated form of the matrix.

\medskip
So up to now  $A$ and $B$ look like
\[ 
A=
\left(\begin{array}{c|c}
\begin{array}{cc}
{\displaystyle 0}  & \!\!{\displaystyle-\E_n}\\[0.5ex]
{\displaystyle\E_n} &{\displaystyle 0}
\end{array} 
& 
\begin{array}{c}
\!\! 0 \\[-0.8ex]\!\! \vdots \\[-0.3ex]\!\! 0
\end{array} \\ \hline
0\, \cdots\,  0 & \!\! 0
\end{array}\right)
\  \ \ 
B=
\left( \begin{array}{c|c|c}
{\displaystyle \widetilde{B}}\! & 
\begin{array}{c}\!\!\! 0\!\!\!\!  \\[-0.8ex] \!\!\!\vdots\!\!\!\!  \end{array} &
\! {\displaystyle -\E_n} \\ \cline{1-1}\cline{3-3}
\multicolumn{3}{c}{\hspace{0.4ex} 0\; \cdots \; 0 \;\cdots\; 0} \\[-0.35ex]\cline{1-1}\cline{3-3}
{\displaystyle \raiseghost{4ex}{0ex}{${\displaystyle \E_n}$}\phantom{-\E_n}}\! & 
\begin{array}{c} \!\!\!\vdots\!\!\!\!  \\[-0.3ex]\!\!\!0\!\!\!\! \end{array} 
&   {\displaystyle 0} 
\end{array}\right).
\]

\noindent\emph{$3^{rd}$ Step:}
We adjust the vectors $(e_0,\ldots,e_{n-1})$ so that $\widetilde{B}=0$ and $A$
stays the same.

\medskip
We note that a transformation of $\CC^{2n+1}$ by 
\[
T=
\left(\begin{array}{c|c}
\begin{array}{cc}
{\displaystyle \E_n}  & {\displaystyle 0}\\[0.5ex]
{\displaystyle t} &{\displaystyle \E_n}
\end{array} 
& 
\begin{array}{c}
\!\! 0 \\[-0.8ex]\!\! \vdots \\[-0.3ex]\!\! 0
\end{array} \\ \hline
0\, \cdots\,  0 & \!\! 1
\end{array}\right)^{\!\!-1}
\ \ \mbox{with}\ t\in \sym(n,\CC) 
\]
does not change $A$ since 
\[
{}^tT^{-1}AT^{-1}={}^tT^{-1}(AT^{-1})=
\left(\begin{array}{c|c}
\begin{array}{cc}
{\displaystyle \E_n}  & {\displaystyle t}\\[0.5ex]
{\displaystyle 0} &{\displaystyle \E_n}
\end{array} 
& 
\begin{array}{c}
\!\! 0 \\[-0.8ex]\!\! \vdots \\[-0.3ex]\!\! 0
\end{array} \\ \hline
0\, \cdots\,  0 & \!\! 1
\end{array}\right)
\!
\left(\begin{array}{c|c}
\begin{array}{cc}
{\displaystyle \!\!-t}  &\!\! {\displaystyle-\E_n}\\[0.5ex]
{\displaystyle\E_n} &{\displaystyle 0}
\end{array} 
& 
\begin{array}{c}
\!\! 0 \\[-0.8ex]\!\! \vdots \\[-0.3ex]\!\! 0
\end{array} \\ \hline
0\, \cdots\,  0 & \!\! 0
\end{array}\right)
=A.
\]

If we denote by $\overline{t}$ resp. $\left|t\right.\in \M(n\times n,\CC)$ the 
matrix that we obtain by deleting the first row resp. column of $t$ 
and adding a row 
of zeroes below resp. a column of zeroes on the right side, we can write down
the transformation of $B$ as follows:
\[\begin{array}{l}
{}^tT^{-1}BT^{-1}={}^tT^{-1}(BT^{-1})=
\left(\begin{array}{c|c}
\begin{array}{cc}
{\displaystyle \E_n}  & {\displaystyle t}\\[0.5ex]
{\displaystyle 0} &{\displaystyle \E_n}
\end{array} 
& 
\begin{array}{c}
\!\! 0 \\[-0.8ex]\!\! \vdots \\[-0.3ex]\!\! 0
\end{array} \\ \hline
0\, \cdots\,  0 & \!\! 1
\end{array}\right)
\left( \begin{array}{c|c|c}
{\displaystyle\! \widetilde{B}-\overline{t}}\! & 
\begin{array}{c}\!\!\! 0\!\!\!\!  \\[-0.8ex] \!\!\!\vdots\!\!\!\!  \end{array} &
 {\displaystyle -\E_n} \\ \cline{1-1}\cline{3-3}
\multicolumn{3}{c}{\hspace{0.65ex} 0\; \cdots \; 0 \;\cdots\; 0} \\[-0.35ex]\cline{1-1}\cline{3-3}
{\displaystyle \raiseghost{4ex}{0ex}{${\displaystyle \E_n}$}\phantom{-\E_n}}\! & 
\begin{array}{c} \!\!\!\vdots\!\!\!\!  \\[-0.3ex]\!\!\!0\!\!\!\! \end{array} 
& {\displaystyle 0} 
\end{array}\right)
\\[4.8ex]
\phantom{{}^tT^{-1}BT^{-1}}=
\left( \begin{array}{c|c|c}
{\displaystyle\! \widetilde{B}-\overline{t}+\left|t\right.}\! & 
\begin{array}{c}\!\!\! 0\!\!\!\!  \\[-0.8ex] \!\!\!\vdots\!\!\!\!  \end{array} &
 {\displaystyle\ \ \ -\E_n \ \ \ } \\ \cline{1-1}\cline{3-3}
\multicolumn{3}{c}{\hspace{0.5ex} 0 \;\ \  \cdots\ \ \; 0 \; \ \ \cdots\ \ \; 0} 
\\[-0.35ex]\cline{1-1}\cline{3-3}
\E_n & 
\begin{array}{c} \!\!\!\vdots\!\!\!\!  \\[-0.3ex]\!\!\!0\!\!\!\! \end{array} 
& {\displaystyle 0} 
\end{array}\right)
.
\end{array}
\]

So to finish this step, we need to show that every antisymmetric matrix 
$\widetilde{B}=(b_{ij})\in\antisym(n,\CC)$ can be written as 
$\overline{t}-\left|t\right.$ for a symmetric matrix $t\in\sym(n,\CC)$. The
entries of $\overline{t}-\left|t\right.$ are 
\[ 
(t_{i+1,j}-t_{i,j+1})_{i j}
\]
where $t_{n+1, i}:=t_{i, n+1}:=0$  for all $i=1\ldots n$. Obviously,
$\overline{t}-\left|t\right. $ is antisymmetric. We set $t_{1 i}:=t_{i1}:=0$
for $i=1\ldots n$ and define recursively for $j$ from $n$ down to $2$
\[
t_{i+1, j}:=t_{i, j+1}+b_{i j} \ \ \ \mbox{for}\ i=1\ldots j-1.
\]
Then by the symmetry of $t$ the whole matrix $t$ is defined and 
$\overline{t}-\left|t\right.=B$
\hfill$\Box$

\bigskip
For the linear system in normal form the rational curve of centers has the 
parametrisation
\[
\begin{array}{l}
c(\lambda\,\colon \mu) 
= 
\ker (\lambda A-\mu B)
= 
\ker\left(
\begin{array}{c@{\,}c@{\,}c|c@{\,}c@{\,}c@{\,}c}
 & & &  -\lambda  &\mu  & &   \\[-0.5ex]
 &\raiseghost{-1ex}{0ex}{\huge 0}& &    &\ddots  &\ddots &   \\[-0.5ex]
 & & &    &  &-\lambda  &\mu   \\ \hline
\lambda & & &    &  & &   \\[-0.5ex]
-\mu &\ddots & &    &  & &   \\[-0.5ex]
 &\ddots &\lambda &    &\raiseghost{2ex}{1ex}{\Huge 0}  & &   \\[-0.5ex]
 & &-\mu &    &  & &   
\end{array} 
\right) \\[11ex]
\phantom{c(\lambda\,\colon \mu)} = (0\,\colon \ldots\,\colon 0\,\colon \mu^n\,\colon \mu^{n-1}\lambda\,\colon \ldots\,\colon \lambda^n).
\end{array}
\]
The $\PP_{2n-2}$ of lines through a center $c_{(\lambda:\mu)}$ is given by
\[
c_{(\lambda:\mu)} \wedge q \ \ \ \ \mbox{where\ }q\in\PP_{2n}\  \mbox{with} \  
{}^t\!c_{(\lambda:\mu)} A q= {}^t\!c_{(\lambda:\mu)} B q=0,
\]
i.e.\ $q$ must be an element of the hyperplane $h_{(\lambda:\mu)}\in\PP_{2n}^*$
\[
\begin{array}{c@{\;=\;}l}
h_{(\lambda:\mu)} & \ker  {}^t\!c_{(\lambda:\mu)} A \cap {}^t\!c_{(\lambda:\mu)} B\\[1.8ex]
 & \ker\left(
\begin{array}{ccccccc}
\mu^n & \mu^{n-1}\lambda &\cdots & \mu \lambda^{n-1} & 0 & \cdots & 0\\
\mu^{n-1}\lambda &\mu^{n-2}\lambda^2 &\cdots & \lambda^n & 0 & \cdots & 0
\end{array}\right)\\[3.5ex]
 & \ker\left(
\begin{array}{ccccccc}
\mu^{n-1} & \mu^{n-2}\lambda &\cdots & \lambda^{n-1} & 0 & \cdots & 0\\
\end{array}\right).
\end{array}
\]
So the hyperplanes $h_{(\lambda:\mu)}$, which are traced out by the $\PP_{2n-2}$ 
of lines through the centers, give rise to a rational normal curve of degree
$n-1$ in the space of hyperplanes containing the center curve. That the hyperplanes
$h_{(\lambda:\mu)}$ contain the center curve could already be seen from the 
Remark \ref{twocenteroncurve}, by which $h_{(\lambda:\mu)}$ must contain any 
line connecting $c_{(\lambda:\mu)}$ with any other point of the center curve.

\bigskip
Now we are ready to study the automorphism group of $\GG(1,2n)\cap L$.

\medskip
Any automorphism $T\in\aut(\GG(1,2n)\cap L)\subseteq\PP\GL(2n+1,\CC)$ has
to map the center curve onto itself 
and also the projective space $P\cong\PP_n$ 
spanned by the center curve onto itself. It is known \cite[10.12]{H}
that the group
of automorphisms of $\PP_n$ fixing a rational normal curve of degree $n$ is 
isomorphic to $\PP\GL(2,\CC)$. If the rational normal curve is given by
\[
\begin{array}{c@{\,}ccc}
c: & \PP_1 & \longrightarrow & \PP_{n}\\
 & (\lambda\,\colon \mu) &\longmapsto & (\mu^n\,\colon \mu^{n-1}\lambda \,\colon  \ldots\,\colon \lambda^n),
\end{array} 
\]
this isomorphism $\PP\GL(2,\CC)\cong\aut(c,\PP_n)$ maps 
\[
t=\left(\begin{array}{cc} a&b\\c&d\end{array} \right)\in\PP\GL(2,\CC)
\]
to $t_{n+1}\in\PP\GL(n+1,\CC)$ where $t_{n+1}$ is the unique matrix such that
\[
t_{n+1}
\left(\begin{array}{c}
\mu^n\\ \mu^{n-1}\lambda \\ \vdots \\ \lambda^n 
\end{array}\right)
=
\left(\begin{array}{c}
(d\mu+c\lambda)^n \\ (d\mu+c\lambda)^{n-1}(b\mu+a\lambda) \\
\vdots \\  (b\mu+a\lambda)^n
\end{array}\right);
\]
for example
\[
t_2=
\left(\begin{array}{cc}
d & c\\ b & a
\end{array} \right)
\qquad
t_3=
\left(\begin{array}{ccc}
d^2&2cd & c^2\\ bd & ad+bc &ac\\b^2 & 2ab & a^2
\end{array} \right).
\]

Applying this to the center curve restricts the form of the transformation $T$ to
\[
T=
\left(\begin{array}{c|c}
* & 0\\ \hline 
{\displaystyle *} & \rule[-1.5ex]{0ex}{4.5ex} {\displaystyle t_{n+1}}
\end{array} \right).
\]

We know further that if $T$ maps $c_{(\lambda:\mu)}$ to 
 $c_{(a\lambda+b\mu :c\lambda+d\mu)}$, then it must map the hyperplane 
$h_{(\lambda:\mu)}$ to  $h_{(a\lambda+b\mu :c\lambda+d\mu)}$. Therefore it
induces also an automorphism on the rational curve $h$ of degree $n-1$ in the dual
projective space $(\PP_{2n}/P)^*$ of hyperplanes containing $P$. Hence $T$ must
be of the form 
\[
T=
\left(\begin{array}{c|c}
\alpha\, \rule[-2.5ex]{0ex}{5.5ex} {}^t\hspace{-0.1ex} t_n^{-1} & 0\\ \hline
* & \rule[-3ex]{0ex}{6.5ex} {\displaystyle \ t_{n+1}\ }
\end{array} \right)
\ \ \mbox{with\ }\alpha\in\CC^*.
\]

We make the following \emph{claim}:
\[
T=
\left(\begin{array}{c|c}
\rule[-1.5ex]{0ex}{4ex}{}^t t_n^{-1} & 0\\ \hline
0 &  \rule[-1.5ex]{0ex}{4.5ex} {\displaystyle t_{n+1}}
\end{array} \right)\in \PP\GL(2n+1,\CC)
\]
is an automorphism of the linear system $\GG(1,2n)\cap L$.

\medskip
\noindent\emph{Proof.} We need to check that for every $t\in\PP\GL(2,\CC)$
\[
{}^tT^{-1}(\lambda A-\mu B) T^{-1}\in \Span\{ A,B\}
\]
for all $\lambda,\mu \in\CC$. 
Since $\PP\GL(2,\CC)$ is a group, this is equivalent to the statement that for 
every  $t\in\PP\GL(2,\CC)$
\[
{}^tT(\lambda A-\mu B) T\in \Span\{ A,B\}
\]
for all $\lambda,\mu \in\CC$. 
Because of the linearity it is enough to do 
this for $(\lambda,\mu)=(1,0)$ and $(0,-1)$. Denoting by $\overline{t_{n+1}}$
resp.  $\underline{t_{n+1}}$ the matrix $t_{n+1}$ with the first resp. last
row deleted, we compute:
\[
\begin{array}{l}
{}^tT(AT)=\!
\left(\!\begin{array}{c|c}
\rule[-1.5ex]{0ex}{4ex} t_n^{-1} & 0\\ \hline
0 &  \rule[-1.5ex]{0ex}{4.5ex} {\displaystyle{}^t t_{n+1}}
\end{array}\! \right)
\!\!%
\left(\!\begin{array}{c|c}
\rule[-1.5ex]{0ex}{5ex} 0 & -\underline{t_{n+1}}\\ \hline
\!\begin{array}{c}{}^tt_n^{-1} \\[-1ex] {\scriptstyle 0 \cdots 0}
\end{array} \!
& \rule[-1.5ex]{0ex}{5.5ex} {\displaystyle 0}
\end{array}\!\right)
\!=\!
\left(\!\begin{array}{c|c}
\rule[-2.5ex]{0ex}{7ex} \mbox{\Large 0} &  -t_n^{-1}\underline{t_{n+1}} \!\\ \hline
{}^t\!(t_n^{-1} \underline{t_{n+1}})\!  &
 \rule[-2.5ex]{0ex}{7.5ex}   \mbox{\LARGE 0}
\end{array}\!\right)
\\[8ex]
{}^tT(BT)=\!
\left(\!\begin{array}{c|c}
\rule[-1.5ex]{0ex}{4ex} t_n^{-1} & 0\\ \hline
0 &  \rule[-1.5ex]{0ex}{4.5ex} {\displaystyle{}^t t_{n+1}}
\end{array}\! \right)
\!\!%
\left(\!\begin{array}{c|c}
\rule[-1.5ex]{0ex}{5ex} 0 & -\overline{t_{n+1}}\\ \hline
\!\begin{array}{c} {\scriptstyle 0 \cdots 0} \\[-0.5ex]{}^tt_n^{-1}
\end{array} \!
& \rule[-1.5ex]{0ex}{5.5ex} {\displaystyle 0}
\end{array}\!\right)
\!=\!
\left(\!\begin{array}{c|c}
\rule[-2.5ex]{0ex}{7ex} \mbox{\Large 0} &  -t_n^{-1}\overline{t_{n+1}}\!\\ \hline
{}^t\!(t_n^{-1} \overline{t_{n+1}})\!  &
 \rule[-2.5ex]{0ex}{7.5ex}   \mbox{\LARGE 0}
\end{array}\!\right)\!.
\end{array}
\]
So, if we show 
\[
\begin{array}{l}
\underline{t_{n+1}}=d\left( t_n {}_0^0 \right)+
c \left( {}_0^0 t_n \right)
\Longrightarrow 
t_n^{-1} \underline{t_{n+1}}=d\left( \E_n {}_0^0 \right)+
c \left( {}_0^0 \E_n \right)
\\[3ex]
\overline{t_{n+1}}=b\left( t_n {}_0^0 \right)+
a \left( {}_0^0 t_n \right)
\Longrightarrow  
t_n^{-1} \overline{t_{n+1}}=b\left( \E_n {}_0^0 \right)+
a \left( {}_0^0 \E_n \right),
\end{array}
\]
where ${}_0^0$ stands for adding a column of zeroes, then
\[
\begin{array}{l}
{}^tT A T=dA+cB\\[1ex]
{}^tT BT=bA+aB.
\end{array}
\]

To show the equality for $\underline{t_{n+1}}$ note that on the one hand
$\underline{t_{n+1}}$ is the unique matrix with
\[
\underline{t_{n+1}}
\left(\begin{array}{c}
\mu^n\\ \mu^{n-1}\lambda \\ \vdots \\ \lambda^n 
\end{array}\right)
=
\left(\begin{array}{c}
(d\mu+c\lambda)^n \\ (d\mu+c\lambda)^{n-1}(b\mu+a\lambda) \\
\vdots \\  (d\mu+c\lambda)(b\mu+a\lambda)^{n-1}
\end{array}\right)
\]
and on the other hand
\[ 
\begin{array}{l}
\left(\begin{array}{c}
  (d\mu+c\lambda)^n \\
  (d\mu+c\lambda)^{n-1}(b\mu+a\lambda) \\
  \vdots \\  
  (d\mu+c\lambda)(b\mu+a\lambda)^{n-1}
\end{array}\right)
=(d\mu+c\lambda)
\left(\begin{array}{c}
  (d\mu+c\lambda)^{n-1} \\ 
  (d\mu+c\lambda)^{n-2}(b\mu+a\lambda) \\
  \vdots \\  
  (b\mu+a\lambda)^{n-1}
\end{array}\right)
\\[7ex]
=
(d\mu+c\lambda)t_n
\left(\begin{array}{c}
  \mu^{n-1}\\ 
  \mu^{n-2}\lambda \\ 
  \vdots \\ 
  \lambda^{n-1} 
\end{array}\right)
=
d t_n
\left(\begin{array}{c}
  \mu^n\\ \mu^{n-1}\lambda \\ \vdots \\ \mu\lambda^{n-1} 
\end{array}\right)
+
c t_n
\left(\begin{array}{c}
  \mu^{n-1}\lambda\\ \mu^{n-2}\lambda^2 \\ \vdots \\ \lambda^n 
\end{array}\right)\\[7ex]
=
\left( 
  d\left(t_n {}_0^0\right)+
  c\left( {}_0^0 t_n\right)
\right)
\left(\begin{array}{c}
   \mu^n\\ \mu^{n-1}\lambda \\ \vdots \\ \lambda^n 
\end{array}\right).
\end{array}
\]

Of course, the proof for $\overline{t_{n+1}}=b\left( t_n {}_0^0 \right)+
a \left( {}_0^0 t_n \right)$ is analogous.
\hfill$\Box$

\bigskip
Given any automorphism of the line system $\GG(1,2n)\cap L$ we can compose
it with one of the above automorphisms such that the composition fixes the
center curve pointwise. So, we can focus our attention to automorphisms of 
the last type.

\begin{lemma}\label{centerfix}
All automorphisms of $\GG(1,2n) \cap L$ that fix the center curve pointwise
are of the form
\[
T=\left(\begin{array}{cc}
\alpha \E_n & 0\\
S & {\displaystyle  \E_{n+1}}
\end{array}\right)
\ \mbox{with\ }\alpha\in\CC^*,\ S\in\M( (n+1)\times n,\CC),
\]
where the matrix $S\in\M( (n+1)\times n,\CC)$ has the same entries along
the minor diagonals, i.e.\ $s_{ij}=s_{kl}$ for $i+j=k+l$.

\medskip
As a group these matrices are isomorphic to the semi direct product
$\CC^{2n}\ltimes \CC^*$,
$\ (s,\alpha)\cdot(s',\alpha')=(\alpha's+s',\alpha\alpha')$.
\end{lemma}
\emph{Proof.} We need only to check the property of $S$ and the group 
structure. $T$ is an automorphism iff
\[
{}^tT^{-1}AT^{-1},{}^tT^{-1}BT^{-1}\in \Span\{A,B\}.
\]

The inverse of $T$ is 
\[
T^{-1}=\left(\begin{array}{cc}
\frac{1}{\alpha} \E_n & 0\\[1ex]
\!- \frac{1}{\alpha} S & {\displaystyle  \E_{n+1}}
\end{array}\right).
\]

Now if $\overline{S}$ resp. $\underline{S}\in\M(n\times n,\CC)$ denote the
matrix $S$ with the first resp. last row deleted, then
\[\begin{array}{l} 
{}^tT^{-1}(AT^{-1})=
\!\!\left(\begin{array}{cc}
\!\frac{1}{\alpha} \E_n\! &\!\!-\frac{1}{\alpha} {}^t\!S \!\\[1ex]
\! 0 & \!\! {\displaystyle\E_{n+1}} \!
\end{array}\right)
\!\!%
\left(\begin{array}{c|c}
\begin{array}{cc}
\!\!\! \frac{1}{\alpha}{\displaystyle  \underline{S}}  &\!\!\! {\displaystyle-\E_n}\!\!\\[1ex]
\!\!\! \frac{1}{\alpha} {\displaystyle \E_n} &{\displaystyle 0}
\end{array}\!\!\!\! 
& 
\begin{array}{c}
\!\! 0 \\[-0.8ex]\!\! \vdots \\[-0.3ex]\!\! 0
\end{array} \\ \hline
0\, \cdots\,  0 & \!\! 0
\end{array}\right)
\\[6.5ex]
\phantom{{}^tT^{-1}(AT^{-1}) }
=
\left(\begin{array}{c|c}
\begin{array}{cc}
\!\!\! \frac{1}{\alpha^2}{\displaystyle(\underline{S}-{}^t\!\underline{S})}  
&\!\! -\frac{1}{\alpha} {\displaystyle\E_n}\\[1ex]
\!\!\frac{1}{\alpha}{\displaystyle\E_n} &{\displaystyle 0}
\end{array}\!\!\!\! 
& 
\begin{array}{c}
\!\! 0 \\[-0.8ex]\!\! \vdots \\[-0.3ex]\!\! 0
\end{array} \\ \hline
0\, \cdots\,  0 & \!\! 0
\end{array}\right)
\\[6.5ex]
{}^tT^{-1}(BT^{-1})=
\left(\begin{array}{cc}
\!\frac{1}{\alpha} \E_n\! & \!\! -\frac{1}{\alpha}  {}^t\!S\!\\[1ex]
\!0 &\!\! {\displaystyle \E_{n+1}}\!\!
\end{array}\right)
\!\!%
\left( \begin{array}{c|c|c}
 \frac{1}{\alpha}{\displaystyle\overline{S}}\! & 
\begin{array}{c}\!\!\! 0\!\!\!\!  \\[-0.8ex] \!\!\!\vdots\!\!\!\!  \end{array} &
\! {\displaystyle -\E_n} \\ \cline{1-1}\cline{3-3}
\multicolumn{3}{c}{\hspace{0.0ex} 0\; \cdots \; 0 \;\cdots\; 0} \\[-0.35ex]\cline{1-1}\cline{3-3}
\frac{1}{\alpha}{\displaystyle  \E_n}\!  & 
\begin{array}{c} \!\!\!\vdots\!\!\!\!  \\[-0.3ex]\!\!\!0\!\!\!\! \end{array} 
&\! {\displaystyle 0} 
\end{array}\right)
\\[8.5ex]
\phantom{{}^tT^{-1}(BT^{-1}) }
=
\left( \begin{array}{c|c|c}
\!\!\frac{1}{\alpha^2}{\displaystyle (\overline{S}-{}^t\overline{S}) }\!\!\! & 
\begin{array}{c}\!\!\! 0\!\!\!\!  \\[-0.5ex] \!\!\!\vdots\!\!\!\!  \end{array} &
\! \ \; -\frac{1}{\alpha} {\displaystyle \E_n}\ \;  \\ \cline{1-1}\cline{3-3}
\multicolumn{3}{c}{\hspace{0.5ex} 0\ \ \;   \cdots \ \ \;  0 \ \ \; \cdots\ \ \; 0} \\[-0.35ex]\cline{1-1}\cline{3-3}
\frac{1}{\alpha}{\displaystyle \E_n}\!  & 
\begin{array}{c} \!\!\!\vdots\!\!\!\!  \\[-0.0ex]\!\!\!0\!\!\!\! \end{array} 
&\! {\displaystyle 0} 
\end{array}\right)\!\!.
\end{array}
\]
Therefore $T$ is an automorphism iff $\underline{S}={}^t\!\underline{S}$ and 
$\overline{S}={}^t\overline{S}$. In other words
\[\left.
\begin{array}{c@{\;=\;}c}
s_{ij} & s_{ji}\\
s_{i+1, j} & s_{j+1, i}
\end{array}
\right.\quad \mbox{for\ }1\le i,j\le n,
\]
so
\[
s_{ij} = s_{ji} = s_{(j-1)+1, i} = s_{i+1, j-1}
\]
for $j>1$ and $i<n$, hence $s_{ij}=s_{kl}$ for $i+j=k+l$.

\medskip
The statement about the group action follows from

\bigskip
\ \hfill
$\displaystyle  
\left(\begin{array}{cc}
\alpha \E_n &\!\! 0\\
S &\!\! \E_{n+1}
\end{array}\right)
\left(\begin{array}{cc}
\alpha' \E_n &\!\! 0\\
S' &\!\! \E_{n+1}
\end{array}\right)
=
\left(\begin{array}{cc}
\!\alpha\alpha' \E_n &\!\! 0\\
\!\alpha'S+S' & \!\!\E_{n+1}
\end{array}\right).
$
\hfill$\Box$

\bigskip
Collecting the results we have

\begin{theorem}\label{autog12nh2}
The automorphism group of $\GG(1,2n)\cap L$ is an extention of $\PP\GL(2,\CC)$
by the semi direct product $\CC^{2n}\ltimes\CC^*$.

\medskip
It is isomorphic to the matrix subgroup of $\,\PP\GL(2n+1,\CC)$ given by
\[
\left(\begin{array}{cc}
\alpha \E_n &\! 0\\
S & \!\E_{n+1}
\end{array}\right)
\left(\begin{array}{cc}
{}^tt_n^{-1} &\! 0\\
0 & \!t_{n+1}
\end{array}\right)
\]
where $\alpha\in\CC^*$, $S\in\M((n+1)\times n,\CC)$ with $s_{i j}=s_{kl}$ for 
$i+j=k+l$ and $t_n\in\aut(h,\PP_{n-1})$ resp. $t_{n+1}\in\aut(c,\PP_n)$ are the 
transformations that are induced by the $\PP\GL(2,\CC)$ action on the rational 
normal curve $h\subset\PP_{n-1}$ resp. $c\subset\PP_n$.
\end{theorem}
\emph{Proof.} It remains to show that the automorphism fixing the center curve pointwise form
a normal subgroup, but that can be easily computed.\hfill$\Box$

\bigskip
\begin{remark}\label{centerdetermine}
An automorphism of $\GG(1,2n)\cap L$ is determined by its action on the lines intersecting
the center curve.
\end{remark}

\medskip
In contrast to that, the line system, i.e.\ the position of the line 
 $L^*\subset \PP(\bigwedge^2\CC^{2n+1})^*$, 
is not determined by these lines, as a simple dimension count shows. Giving these lines is 
equivalent to giving the two rational curves $c\subset\PP_{2n}$ and 
$h\subset \PP_{2n}/P\cong\PP_{n-1}$ and a correspondence between them, so that 
we have the following dimension count
\[
(2(2n+1)-4)+(2n-4) +3 < \dim \GG(1,\PP({\textstyle \bigwedge}^2\CC^{2n+1}))= 
2\left({ 2n+1 \choose 2}-2\right).
\]
\emph{Proof of the remark.} We need to show that only the identity fixes 
these lines one by one.
First a transformation $T$ that fixes the lines must fix the center curve, 
hence by the Lemma
\ref{centerfix} it is of the form
\[
T=
\left(\begin{array}{cc}
\alpha \E_n &\! 0\\
S & \!\E_{n+1}
\end{array}\right).
\]

We compute the induced action $\widetilde{T}$ of $T$ on
$\{ l\in\GG(1,2n)\cap L\mid c_{(0:1)} \in l\}$, 
the $\PP_{2n-2}$ of lines through $c_{(0:1)}=e_n$.
A line $l\in\GG(1,2n)$ through $c_{(0:1)}$ will be 
in the line system $\GG(1,2n)\cap L$ iff
it lies in the hyperplane $h_{(0:1)}=\ker(1\,\colon 0\,\colon \ldots\,\colon 0)$. Therefore the $\PP_{2n-2}$ of lines
through $e_n$ is given by
\[
e_n\wedge x \ \quad \mbox{with\ }x\in\PP(\Span\{e_1,\ldots,e_{n-1},e_{n+1},\ldots,e_{2n}\}).
\]
Using $(e_1\wedge e_n,\ldots,e_{n-1}\wedge e_n,e_{n+1}\wedge e_n,\ldots,e_{2n}\wedge e_n)$ as a basis, the induced action 
$\widetilde{T}$ is
\[
\widetilde{T}=
\left(\begin{array}{cc}
\alpha \E_{n-1} &\! 0\\[0.5ex]
|\overline{S} & \!\E_n
\end{array}\right).
\]
Here $|\overline{S}$ denotes the matrix $S$ with the first row and column deleted.
In order to have $\widetilde{T}=\E_{2n-1}$, we must have $\alpha=1$ and 
$|\overline{S}=0$.

\medskip
The same computation for the lines through $c_{(1:0)}=e_{2n}$ 
yields $\alpha=1$ and 
$\underline{S}|=0$ from which $S=0$ and the remark follow.
\hfill $\Box$

\bigskip
For the rest of the section we analyze the action of the automorphism 
group on the line system $\GG(1,2n)\cap L$.
We start with  $\GG(1,4)\cap L$.

\begin{proposition}
The action of $\aut( \GG(1,4)\cap L)$ on the lines has four orbits:
\begin{enumerate}
\item tangents of the center conic
\item secants of the center conic
\item lines through the center conic that do not lie in the plane of the
center conic
\item lines that do not intersect the plane of the center curve.
\end{enumerate}
\end{proposition}
\emph{Proof.} Since any automorphism maps the center conic onto itself, it 
is clear by the geometric description that all the mentioned lines lie in 
different orbits.

\medskip
Any line in the plane $P$ of the center conic intersects the conic twice,
so by the Remark \ref{twocenteroncurve} it is a member of the line system. 
Since the
automorphism group acts like $\aut(c,P)\cong\PP\GL(2,\CC)$ on the plane $P$,
the first two orbits are obvious.

\medskip
To see that the lines of 3) form one orbit, we have to exhibit an automorphism
that given two lines of 3) maps one onto the other. Since the $\PP\GL(2,\CC)$
part of the automorphism group acts transitively on the center conic, we may 
assume that both lines pass through the same point of the center conic, say
$e_2=c_{(0:1)}$. Now the induced action $\widetilde{T}$ of an automorphism 
$T$ fixing the center conic pointwise on the $\PP_2$ of lines through
$e_2$ was computed in the proof of the Remark \ref{centerdetermine} as
\[
\widetilde{T}=
\left(\begin{array}{ccc}
\alpha & 0 & 0 \\
f & 1 & 0 \\
g & 0 & 1
\end{array} \right)
\ \ \mbox{with\ }\alpha\in\CC^*;\ f,g\in \CC.
\]
These transformations act transitively on 
$\PP_2\setminus\PP(\Span\{\widetilde{e_1},\widetilde{e_2}\})$,
where the line $\PP(\Span\{\widetilde{e_1},\widetilde{e_2}\})$ corresponds
to the lines through $e_2$ that lie in the plane of the center conic.

\medskip
The lines of 4) are all the remaining lines since there are no lines 
that intersect the plane $P$ of the center conic but not the conic $c$
itself. This is clear because the $\PP_1$ of lines through a point 
$p\in P\setminus c$ is formed by the lines through $p$ in the plane $P$,
so there can be no other line.

\medskip
Finally, we have to  show that the lines of 4) form one orbit.
By a small computation one checks that 
$\aut(\GG(1,4)\cap L)\subset \PP\GL(5,\CC)$ acts transitively 
on $\PP_4\setminus P$. So, it suffices to show that the line $e_0\wedge e_1$ 
can be mapped to any other line through $e_0$  by an automorphism. Any of 
these lines can be written as
\[
e_0\wedge (e_1+\beta e_4)
\quad \ \mbox{with\ }\beta \in \CC,
\]
and the automorphism 
\[
T=
\left( \begin{array}{cc|ccc}
1  & 0 & 0 &0 &0 \\
0 & 1 &0 & 0 &0 \\ \hline
0 & 0 & 1 & 0 &0 \\
0 & 0 &0 & 1 &0 \\
0 & \beta &0 & 0 &1
\end{array}\right) 
\]
will take $e_0\wedge e_1$ to it.
\hfill $\Box$

\begin{proposition}
The automorphism group  acts quasihomogeneously on $\GG(1,6)\cap H^2$.
\end{proposition}
\emph{Proof.} For this it is enough to show that the stabiliser of the line 
$l=e_0\wedge e_2$ is a 2-dimensional subgroup since then
\[
\begin{array}{c@{\,=\,}l}
 \dim \mathrm{Orbit}(l) & \dim \aut (\GG(1,6)\cap H^2) -\dim \mathrm{Stab}(l)
=10-2=8\\
 & \dim \GG(1,6)\cap H^2.
\end{array}
\]

If we normalize the $t\in\PP\GL(2,\CC)$ by $\det t=1$, every 
$T\in \aut (\GG(1,6)\cap H^2)$  can be written by the Theorem \ref{autog12nh2}
as
\[
\begin{array}{l}
T=
\left( \begin{array}{ccc|cccc}
\alpha &  &0  &     &  &  & \\
 & \alpha &   &     & \raiseghost{3.5ex}{-1ex}{\Huge 0} &  &  \\
0 &  &\alpha  &     &  &  &  \\ \hline
e & f & g     &    1&  & \raiseghost{3ex}{-2ex}{\LARGE 0} &  \\
f & g & h     &     & 1&  &  \\
g & h & i     &     &  & 1&  \\
h & i & j     &  \raiseghost{3.5ex}{0.5ex}{\LARGE 0}  &  &  & 1 
\end{array} \right)
\left(\begin{array}{cc}
{}^tt_3^{-1} &\! 0\\
0 & \!t_4
\end{array}\right)\\[11ex]
\ \ \mbox{with}\ {}^tt_3^{-1}=
\left( \begin{array}{ccc}
a^2 & -ab &b^2   \\
-2ac& ab+cd &-2bd \\
c^2 & -cd &d^2  
\end{array} \right).
\end{array}
\]

To compute the stabilizer we start by looking only at the first three 
entries of 
\[
\begin{array}{l}
Te_0=(\alpha a^2 , -2\alpha ac, \alpha c^2,\ldots)\\[1ex]
Te_2=(\alpha b^2 , -2\alpha bd, \alpha d^2,\ldots).
\end{array}
\]
Since we must have $Te_0,Te_2\in l$, $ac=bd=0$ follows.
By $\det t=ad-bc=1$ we have the two possibilities $b=c=0$, $d=a^{-1}$ and
$a=d=0$, $c=-b^{-1}$. We examine only the first case, the second being similar.
Now we have
\[
\begin{array}{l}
Te_0=a^2(\alpha , 0, 0 ,e,f,g,h)\\[1ex]
Te_2=a^2(0,0,\alpha,g,h,i,j).
\end{array}
\]
From $Te_0,Te_2\in l$ we conclude $e=f=g=h=0$ resp. $g=h=i=j=0$.
Therefore, including the case $(a=d=0,\ c=-b^{-1})$, the stabilizer of $l$ is 

\bigskip
\hfill
$%\displaystyle
\mathrm{Stab}(l)=\left\{
\left( \begin{array}{c@{\,}c@{\,}c@{\,}c@{\,}c@{\,}c@{\,}c}
\alpha a^2\! & & & & & & \\
 &0  & & & &\raiseghost{-1ex}{-1ex}{\Huge 0} & \\
& & \alpha a^{-2}\! & & & &  \\
 & & &a^{-3} & & & \\
 & & & &\!\!a^{-1}\, & & \\
 &\raiseghost{-1ex}{0ex}{\Huge 0} & & & &\!\!\!\!\!a & \\
 & & & & & &a^3
\end{array} \right),
\left( \begin{array}{c@{\,}c@{\,}c@{\,}c@{\,}c@{\,}c@{\,}c}
\!\!0 & &\!\alpha b^2 & & & & \\
 &0 & & & &\raiseghost{0ex}{-0.5ex}{\LARGE 0} & \\
\!\!\alpha b^{-2} & &0 & & & & \\
 & & & & & &\!\!\!\!\!-b^{-3} \\
 & & & &  &\;b^{-1}\! & \\
 &\raiseghost{0ex}{1ex}{\LARGE 0} & & &\!\!\!-b & & \\
 & & &\!b^3 &   & &\raiseghost{0ex}{1ex}{\Large 0}
\end{array} \right)\right\}.$
\hfill
$\Box$

\bigskip

\begin{proposition}
For $n\ge4$ the action of  the automorphism group 
on  $\GG(1,2n)\cap H^2$ is not quasihomogeneous.
\end{proposition}
\emph{Proof.} We  project the $\PP_{2n}$ 
from the space $P$ of the center curve onto $\PP_{2n}/P\cong \PP_{n-1}$.
This projects the lines of $\GG(1,2n)\cap H^2$ not intersecting $P$ 
surjectively onto the lines $\GG(1,\PP_{2n}/P)$ of $\PP_{2n}/P$. The
automorphisms of  $\GG(1,2n)\cap H^2$ induce automorphisms of
$\PP_{2n}/P$. As matrices these are the upper left $n\times n$ matrices of the 
matrices of Theorem \ref{autog12nh2}, i.e.\ they are of the form
${}^tt_n^{-1}$. So, as a group this induced automorphism group is isomorphic to
$\aut(h^*,\PP_{2n}/P)\cong \PP\GL(2,\CC)$. If 
 $\aut(\GG(1,2n)\cap H^2)$ acts quasihomogeneously, then 
$\aut(h^*,\PP_{2n}/P)$ would have to act quasihomogeneously on 
$\GG(1,\PP_{2n}/P)\cong\GG(1,n-1)$, but this contradicts

\bigskip
\ \hfill
$\displaystyle
\dim \PP\GL(2,\CC)=3<\dim \GG(1,n-1)=2n-4.$
\hfill $\Box$

\section{$\GG(1,4) \cap H^3$}

Let $L=H^3$ be a general 3-codimensional subspace of 
$\PP(\bigwedge\nolimits^2 \CC^5) \cong \PP_9$. 
To $L$ corresponds the plane 
$L^*=\PP(\lambda A+\mu B+\nu C) \subseteq \PP(\bigwedge\nolimits^2\CC^5)^*$
of hyperplanes containing $L$. Since the locus of antisymmetric matrices of  
corank 3 is 3-codimensional in $\PP(\bigwedge\nolimits^2 \CC^5)^*$ by 
Corollary \ref{dualG}, $L^*$ does not contain any.
Therefore to each of the hyperplanes $H_{(\lambda:\mu:\nu)}\subset L$
corresponds a unique center $c_{(\lambda:\mu:\nu)} \in \PP_4$. 
In complete analogy to the last case we get

\begin{lemma} \label{centerofg14}
The centers  $c_{(\lambda:\mu:\nu)}$ are those points of $\PP_4$ through
which there passes a $\PP_1$  of lines of the line system $\GG(1,4) \cap L$. Through all 
the other points passes a unique line.
\end{lemma}

\begin{proposition}\label{centervero} The map of centers 
\[ \begin{array}{rccc}
c: & L^* \cong \PP_2 & \longrightarrow & \PP_4\\[0.5ex]
 & (\lambda:\mu:\nu) & \longmapsto & c_{(\lambda:\mu:\nu)}=\ker (\lambda A+\mu B+\nu C)
\end{array} \]
is an embedding of $\PP_2$ in $\PP_4$ of degree 2, i.e. its image is a smooth projected 
Veronese surface.
\end{proposition}

\begin{remark} Any line that contains three centers is in the line system.
\end{remark}
\emph{Proof.} Let $c(p_0)$, $c(p_1)$ and $c(p_2)$ with 
$p_0,p_1,p_2 \in L^*$ be the three centers on the line $l$. 
By the definition of the centers we have $l\in H_{p_i}$.
Since $c$ maps lines in $L^*$ onto conics in $\PP_4$, the three 
points   $p_0,p_1,p_2$ do not lie on a line, hence they span $L^*$. So 
$ l\in H_{p_0} \cap H_{p_1} \cap H_{p_2} =L$. \hfill $\Box$

\bigskip
From the statements we get a complete picture of the lines of 
$\GG(1,4) \cap H^3$ in $\PP_4$. 
We define the trisecant variety $\Tri(X)$ of a variety 
$X\subseteq \PP_N$ by:
\[
\Tri(X):=\overline{\{ l\in\GG(1,N) \mid \#(X \cap l)\ge 3 \}}
\]
Then we have

\begin{corollary}[Castelnuovo] $\GG(1,4) \cap L$ is  the trisecant variety  of 
the smooth projected Veronese surface $\im c \subset \PP_4$.
\end{corollary}
\emph{Proof.}(see \cite{C} or \cite[X, 4.4]{SR})
By the remark above the trisecant variety is contained in the irreducible 
variety  $\GG(1,4) \cap L$. So it is enough to show that both varieties 
have the same dimension. The Lemma  \ref{centerofg14} together with the 
Proposition \ref{centervero} shows that there 
is an unique line of $\GG(1,4) \cap L$ through a general point of $\PP_4$ 
and that  $\GG(1,4) \cap L$ is the closure
of such lines. The same statement for the trisecant variety is classical 
 \cite[VII,3.2]{SR}. Hence both varieties have dimension three. \hfill $\Box$

\bigskip
The general trisecant intersects the projected Veronese 
surface $\im c$ in three different 
points. Their inverse image under $c$ are triples of 
points in $\PP_2$ that have to fulfill
some conditions since there is only a 3-dimensional family of 
these triples. To see what 
these conditions are, we recall some facts about the Veronese surface \cite{H}.

\medskip
The Veronese surface $V$ is the image of the embedding 
\[ \begin{array}{rccc}
\upsilon :& \PP_2=\PP(\CC^3) & \longrightarrow & \PP(\sym^2 \CC^3)\\[0.5ex]
 & \PP(v) & \longmapsto & \PP( v\cdot v)\, .
\end{array} \]

Its secant variety consists of the points of $\PP(\sym^2 \CC^3)$ that are the product of
two vectors of $\CC^3$,
\[ \sec(V)=\left\{ \PP(v \cdot w) \mid v,w\in \CC^3\setminus \{0\}\right\}.\]

The projected Veronese surface will be smooth -- like in our case -- iff the center of 
projection $P$ is not in the secant variety.

\bigskip
An intersection of the Veronese surface $V$ with a hyperplane $H\in \PP(\sym^2 \CC^3)^*$ 
gives the conic $\upsilon^{-1}(V\cap H)\subset \PP_2$ which is described by the equation $H$ 
if we identify $\PP(\sym^2 \CC^3)^*$ with the polynomials of degree 2 modulo $\CC^*$.
The conics that we get as hyperplane sections of the projected Veronese surface are 
precisely the conics that we get as hyperplane sections of the Veronese surface by 
hyperplanes that contain the projection center $P$. 
So these conics fulfill one linear condition given by $P$.
We view  $P\in \PP(\sym^2 \CC^3)$ as a conic $C^*_P$ in $\PP_2^*$.
Since $P$ does not lie in the secant variety of the Veronese surface, $P$ is not the product
of two elements of $\CC^3$, hence $C^*_P$ is not the union of two lines. 
Therefore it is smooth. 
We denote the dual conic of  $C^*_P$ by $C_P\subset \PP_2$.

\medskip
Now, three different points of the projected Veronese surface $\im c\subset \PP_4$ 
lie on a line, 
the  trisecant, iff any hyperplane that contains two of them contains all three. 
Under the inverse of the embedding $c$ that means the following on the $\PP_2$:

\medskip
Three different points of $\PP_2$ are the inverse image $c^{-1}(l)$ of a trisecant $l$ 
of the projected Veronese surface $\im c$ iff all conics that fulfill the linear condition
given by $P$ (or equivalently by $C_P$) and pass through two of the points pass through
all three of them.

\medskip
The propositions in the appendix tell us that these triples of points are the vertices of the
non-degenerated polar triangles of the conic $C_P$. By a continuity argument the trisecants
that are tangent to $\im c$ at one point and intersect it in another correspond to the 
degenerated polar triangles, and the trisecants that intersect $\im c$ in only one point 
``with multiplicity three''  correspond to a triple point on the conic $C_P$.
We also see that there are no 4-secants. Since if there is one, there would be four points 
in $\PP_2$ such that any three of them build a different polar triangle. But this is 
impossible because a polar triangle is already determined by two of its vertices.

\bigskip
With this geometric description it is easy to compute the automorphism group of 
$\GG(1,4)\cap H^3$. Any automorphism as a projective linear transformation of $\PP_4$
maps by definition the trisecants of the projected Veronese surface $\im c$ onto themselves.
Further, it must fix the projected Veronese surface $\im c$, since $\im c$ is the union of
the centers. In fact, the automorphism is already determined by its action on $\im c$ 
since this action on $\im c$ determines the images on the trisecants.
Under the inverse of the embedding $c$, this automorphism of $\im c$ 
preserving the 
trisecants corresponds to an automorphism of $\PP_2$ preserving the polar 
triangles of the conic $C_P \subset \PP_2$. Such an automorphism of 
$\PP_2$ maps the 
degenerated polar triangles onto themselves. In particular, it maps the tangents to the 
conic $C_P$ onto themselves. Therefore it has to fix the conic $C_P$.

\medskip
So we have seen how an automorphism of $\GG(1,4)\cap L$ induces an unique automorphism
of $\PP_2$ fixing $C_P$, hence an automorphism of $C_P$ since 
$\aut(C_P,\PP_2)\cong \aut(C_P)\cong \PP \GL(2,\CC)$.

\medskip
On the other hand, any projective linear transformation of $\PP_2$ 
that fixes the conic $C_P$
preserves the polar triangles of $C_P$. Therefore it induces via the 
embedding $c$ an automorphism of the projected Veronese surface $\im c$ 
that preserves triples of points
that lie on a line. So it defines an automorphism of the trisecants of $\im c$, which is 
the same as an automorphism of $\GG(1,4)\cap L$.

\medskip
We summarize:

\begin{theorem} The automorphism group of $\GG(1,4)\cap H^3$ is isomorphic 
to $\PP \GL (2,\CC)$.
\end{theorem}

The description of the orbits of this automorphism group follows immediately.

\begin{proposition} The action of $\aut ( \GG(1,4)\cap H^3 )$ on the linear system 
$\GG(1,4)\cap H^3$ has three orbits:
\begin{enumerate}
\item trisecants of the projected Veronese surface that intersect it in three points
\item trisecants that are tangent to the projected Veronese surface at one point and
intersect it in another
\item trisecants that intersect the projected Veronese surface in only one point ``with 
multiplicity three''.
\end{enumerate}
\end{proposition}
\emph{Proof.} By what was said above, 
this is equivalent to the classical statement that
the action of group $\aut (C_P,\PP_2)$ on the polar triangles 
has three orbits: 
the non-degenerated triangles, the degenerated ones and the triple 
points on $C_P$.\hfill$\Box$

\section{Appendix: Polar Triangles}
Here we prove the needed propositions about polar triangles. The whole appendix may be seen as a modern exposition of \cite[348]{SF}.
First we recall the basic definitions. 

\medskip
Let $C_A$ be a smooth conic in $\PP_2$, which is given by the quadratic 
equation ${}^t\!x A x=0$, where $A\in \GL(3,\CC)$ is  a symmetric, invertible 
matrix. Then $C_A$ induces a polarity $P$ by
\[\begin{array}{rccc}
P:&\PP_2 & \longrightarrow & \PP_2^* \\[0.5ex]
 & \PP(x) & \longmapsto & \PP({}^t\!x A)\, .
\end{array} \] 

For a point $p\in \PP_2$ the line $P(p)$ is called the polar of $p$ and 
$p$ the pole of $P(p)$.

\medskip
A polar triangle is given by three points, at least two of which are 
different, such that the polar of each point contains the other two points,
i.e. $(p,q,r)$ is a polar triangle if ${}^t\!p A q={}^t\!q A r={}^t\!r A p=0$.
The sides of the triangle are the polars of the points.
In the  non-degenerated case when all three points are different, 
the three points cannot lie on a line and therefore span the whole $\PP_2$.
In the degenerated case, $(p,p,q)$, $p$ lies on the conic and $q$ on the 
tangent to the conic at the point $p$. The sides are the polar of $q$ and 
twice the tangent.

\bigskip
%\ 
%{}\hfill 
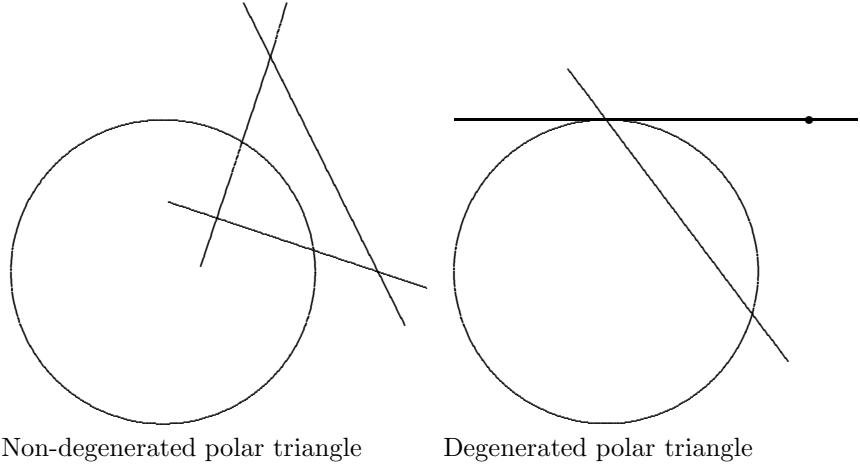
\parbox[b]{5.8cm}{
\mbox{
\font\thinlinefont=cmr5
\begingroup\makeatletter\ifx\SetFigFont\undefined%
\gdef\SetFigFont#1#2#3#4#5{%
  \reset@font\fontsize{#1}{#2pt}%
  \fontfamily{#3}\fontseries{#4}\fontshape{#5}%
  \selectfont}%
\fi\endgroup%
\mbox{\beginpicture
\setcoordinatesystem units <0.750000cm,0.7500000cm>
\unitlength=0.750000cm
\linethickness=1pt
\setplotsymbol ({\makebox(0,0)[l]{\tencirc\symbol{'160}}})
\setshadesymbol ({\thinlinefont .})
\setlinear
%
% Fig ELLIPSE
%
\linethickness= 0.500pt
\setplotsymbol ({\thinlinefont .})
\ellipticalarc axes ratio  2.697:2.697  360 degrees 
        from  5.554 21.907 center at  2.857 21.907
%
% Fig POLYLINE object
%
\linethickness= 0.500pt
\setplotsymbol ({\thinlinefont .})
\plot  3.810 22.860  4.763 25.718 /
%
% Fig POLYLINE object
%
\linethickness= 0.500pt
\setplotsymbol ({\thinlinefont .})
\plot  3.810 22.860  6.668 21.907 /
%
% Fig POLYLINE object
%
\linethickness= 0.500pt
\setplotsymbol ({\thinlinefont .})
\plot  4.763 25.718  6.668 21.907 /
%
% Fig POLYLINE object
%
\linethickness= 0.500pt
\setplotsymbol ({\thinlinefont .})
\plot  6.668 21.907  7.144 20.955 /
%
% Fig POLYLINE object
%
\linethickness= 0.500pt
\setplotsymbol ({\thinlinefont .})
\plot  4.286 26.670  4.763 25.718 /
%
% Fig POLYLINE object
%
\linethickness= 0.500pt
\setplotsymbol ({\thinlinefont .})
\plot  4.763 25.718  5.048 26.670 /
%
% Fig POLYLINE object
%
\linethickness= 0.500pt
\setplotsymbol ({\thinlinefont .})
\plot  6.668 21.907  7.525 21.622 /
%
% Fig POLYLINE object
%
\linethickness= 0.500pt
\setplotsymbol ({\thinlinefont .})
\plot  3.810 22.860  2.953 23.146 /
%
% Fig POLYLINE object
%
\linethickness= 0.500pt
\setplotsymbol ({\thinlinefont .})
\plot  3.810 22.860  3.524 22.003 /
\linethickness=0pt
\putrectangle corners at  0.144 26.695 and  7.550 19.196
\endpicture}
}\\
Non-degenerated polar triangle}
\hfill 
\parbox[b]{6.38cm}{
\mbox{
\font\thinlinefont=cmr5
\begingroup\makeatletter\ifx\SetFigFont\undefined%
\gdef\SetFigFont#1#2#3#4#5{%
  \reset@font\fontsize{#1}{#2pt}%
  \fontfamily{#3}\fontseries{#4}\fontshape{#5}%
  \selectfont}%
\fi\endgroup%
\mbox{\beginpicture
\setcoordinatesystem units <1.06066cm,1.06066cm>
\unitlength=1.06066cm
\linethickness=1pt
\setplotsymbol ({\makebox(0,0)[l]{\tencirc\symbol{'160}}})
\setshadesymbol ({\thinlinefont .})
\setlinear
%
% Fig ELLIPSE
%
\linethickness= 0.500pt
\setplotsymbol ({\thinlinefont .})
\ellipticalarc axes ratio  1.905:1.905  360 degrees 
        from  4.763 22.860 center at  2.857 22.860
%
% Fig ELLIPSE
%
\linethickness= 0.500pt
\setplotsymbol ({\thinlinefont .})
\put{\makebox(0,0)[l]{\circle*{ 0.1}}} at  5.400 24.765
%
% Fig POLYLINE object
%
\linethickness= 0.500pt
\setplotsymbol ({\thinlinefont .})
\putrule from  0.953 24.765 to  6.0 24.765
%\putrule from  0.953 24.765 to  6.668 24.765
% Fig POLYLINE object
%
\linethickness= 0.500pt
\setplotsymbol ({\thinlinefont .})
\plot  2.857 24.765  5.133 21.738 /
%
% Fig POLYLINE object
%
\linethickness= 0.500pt
\setplotsymbol ({\thinlinefont .})
\plot  2.857 24.765  2.381 25.400 /
\linethickness=0pt
\putrectangle corners at  0.927 25.425 and  6.693 20.940
\endpicture}
}\\
Degenerated polar triangle}
%\hfill
{}%\ 

\begin{proposition} Let $C_A=\{x\in \PP_2 \mid {}^t\!x A x=0\}$ be a 
smooth conic and 
$C_A^*=\{x\in \PP_2^* \mid {}^t\!x A^{-1} x=0\}$ its dual conic. 
Further, let 
$C_B=\{x\in \PP_2 \mid {}^t\!x B x=0\}$ be a conic such that 

\bigskip
\ \hfill $\displaystyle \sum_{i,j=0}^2 a^{ij} b_{ij}=0$\hfill $(*)$

\bigskip
\noindent where $A^{-1}=(a^{ij}),B=(b_{ij})\in \sym(3,\CC)$. Finally, let $(p,q,r)$ be a 
polar triangle of $C_A$ then:

\medskip
If two of the three points $p,q,r$ lie on the conic $C_B$, then also the third.
\end{proposition}

In the case of a degenerate polar triangle, $(p,p,q)$, 
the condition that $C_B$ contains $p$ 
twice means that $C_B$ contains $p$ and $C_B$ is 
either singular in $p$ or its tangent in 
$b$ is the polar of $q$.

\medskip
\noindent\emph{Proof.} One can show that all the properties in the statement of the proposition are 
independent of the choice of coordinates, so we may pick nice ones.
We have to distinguish between the two cases of the polar triangle being degenerated or not. 
We treat the case of the non-degenerated polar triangle first. 

\medskip
By a suitable choice of 
coordinates we may assume that $p=(1\,\colon 0\,\colon 0)$, 
$q=(0\,\colon 1\,\colon 0)$ and $r=(0\,\colon 0\,\colon 1)$. Then the assumption
that $(p,q,r)$ is a polar triangle of $C_A$ translates into
\[ 
{}^t\!p A q={}^t\!q A r={}^t\!r A p=0 
\Longleftrightarrow a_{0 1}=a_{1 2}=a_{0 2}=0 . 
\]

By a scaling of the coordinates we can achieve that $ a_{0 0}=a_{1 1}=a_{2 2}=1$, so
that $C_A=\{x\in \PP_2 \mid x_0^2+x_1^2+x_2^2=0\}$. Then the condition 
$(*)$ reads $ b_{0 0}+b_{1 1}+b_{2 2}=0$. If the two points $p$ and 
$q$ are on the conic $C_B$, we 
have ${}^t\!p B p =b_{0 0}=0$ and  ${}^t\!q B q =b_{1 1}=0$. By $(*)$ we see 
$0=b_{2 2}= {}^t\!r B r$, i.e.\ the third point $r$ lies also on the conic $C_B$.

\medskip
Now we treat the case of the degenerated polar triangle $(p,p,q)$. 
We choose coordinates such
that $C_A=\{ x\in \PP_2 \mid  x_0^2+x_1^2+x_2^2=0\}$ and $p=(1\,\colon i\,\colon 0)$. The point $q\neq p$ 
must lie on the tangent to $C_A$. So it has coordinates $q=(\lambda\,\colon i\lambda\,\colon 1)$, and its polar is spanned by $p$ and $(1\,\colon 0\,\colon -\lambda)$. Now using the assumptions

\bigskip
\hspace{3cm}$\displaystyle  b_{0 0}+b_{1 1}+b_{2 2}=0$ \hfill $(*)$

\medskip
\hspace{3cm}$\displaystyle p\in C_B \Longleftrightarrow {}^t\!p B p=0 
\Longleftrightarrow  b_{0 0}+2 i b_{0 1} -b_{1 1}=0,  $ \hfill $(**)$

\bigskip
\noindent we have to show 
\[ 
q\in B \Longleftrightarrow C_B \mbox{\ singular\ in\ }p \ \ \mbox{\ or\ } \ \ 
\TT_pC_B=\; \mbox{polar\ of\ }q.
\]

We rewrite this as 
\[ {}^t\!q B q= {\vphantom{
\left({\textstyle \begin{array}{c}\lambda\\ i\lambda \\1 \end{array}}\right)}
}^t\!\!\!\! 
\left({\textstyle \begin{array}{c}\lambda\\ i\lambda \\1 \end{array}}\right) 
\! B \!
\left({\textstyle \begin{array}{c}\lambda\\ i\lambda \\1 \end{array}}\right) =0
\Longleftrightarrow 
{}^t\! p B \! 
\left({\textstyle \begin{array}{c}1 \\ 0 \\ \!\! -\lambda  \end{array}}\right)=
 {\vphantom{
\left(\begin{array}{c}1\\ i \\1 \end{array}\right)} 
}^t\!\!\!\! 
\left(\begin{array}{c}1\\ i \\0 \end{array}\right) 
\! B\!
\left(\begin{array}{c}1 \\ 0 \\\!\!-\lambda \end{array}\right)=0. 
\]
But this is true since $-2$ times the left hand side plus $(\lambda^2+1)$ times $(**)$
plus $(*)$ gives the right hand side. \hfill $\Box$

\bigskip
Now we prove the converse of the last proposition.

\begin{proposition} Given a smooth conic $C_A=\{x\in \PP_2 \mid {}^t\!x A x=0\}$ 
\[{\cal B}:= \left\{ C_B=\{x\in \PP_2 \mid {}^t\!x B x=0\} \left|  
\sum_{i,j=0}^2 a^{ij} b_{ij}=0\right.\right\} \]
is a four dimensional family of conics. Let $p,q,r\in \PP_2$ be three points, at least two of 
which are different, with the property that if two of them lie on a conic $C_B\in {\cal B}$ then also the third.

\medskip
Then $(p,q,r)$ is a polar triangle of $C_A$.
\end{proposition}
\emph{Proof.} We will show that if $(p,q,r)$ is not a polar triangle then there exits a 
$C_A\in {\cal B}$ for which this property is violated. We have to treat several cases.

\medskip
First let the three points be all different, then they cannot lie on a line. Because if they
would, the conics in the at least two dimensional family
\[ {\cal B}_{p,q}:=\{ C_B \in {\cal B}\mid  p,q\in C_B\} \]
of conics of ${\cal B}$ passing through $p$ and $q$ must split off the line through the
three points. If we pick coordinates such that this line is given by $\{x_2=0\}$, then
${\cal B}_{p,q}$ must be 
\[ {\cal B}_{p,q} =\left\{ \mathrm{V}( x_2(\lambda_0 x_0+\lambda_1 x_1+\lambda_2 x_2)) \mid
(\lambda_0\,\colon \lambda_1\,\colon \lambda_2)\in \PP_2 \right\}.\]
This means that the conics $C_B$ with the matrices
\[ B=\left({\textstyle \begin{array}{ccc}
0 & 0 &  b_{0 2}\\
0 & 0 &  b_{1 2}\\
b_{0 2} & b_{1 2} & b_{2 2}
\end{array}}\right) 
\quad \mbox{with} \ b_{0 2}, b_{1 2}, b_{2 2}\in \CC \]   
are all in ${\cal B}$. Hence the matrix $A^{-1}$ must be of the type
\[ A^{-1}=\left({\textstyle \begin{array}{ccc}
a^{00} & a^{01}   &  0 \\
a^{01} & a^{11}  &  0 \\
0 & 0 & 0
\end{array}}\right) 
\quad \mbox{with} \ a^{0 0}, a^{0 1}, a^{1 1}\in \CC, \]  
but this contradicts the invertibility of $A^{-1}$.

\medskip
Now since $p,q,r$ span the $\PP_2$ we may pick coordinates such that 
$p=(1\,\colon 0\,\colon 0)$, $q=(0\,\colon 1\,\colon 0)$ and $r=(0\,\colon 0\,\colon 1)$. That $(p,q,r)$ is not a polar triangle of $C_A$ 
means that ${}^t\!p A q=a_{01}\not= 0$, ${}^t\!q A r=a_{12}\not= 0$ or 
${}^t\!r A p=a_{02}\not= 0$. 
Assuming $\det A=1$ we conclude that not all of the 
$a^{0 1}=a_{0 2}a_{1 2}-a_{0 1}a_{2 2}$,
$a^{0 2}=a_{0 1}a_{1 2}-a_{0 2}a_{1 1}$ and
$a^{1 2}=a_{0 1}a_{0 2}-a_{0 0}a_{1 2}$
can be zero. If for example $a^{0 2}\not= 0$, then

\[ B=\left({\textstyle \begin{array}{ccc}
0 & 0 & - a^{1 1}\\
0 & 2 a^{0 2} & 0\\
- a^{1 1} & 0 & 0
\end{array}}\right)\]
gives a conic $C_B\in {\cal B}$ that contains the points $p$ and $q$, but not $r$.

\medskip
Now let us look at the case where two of the points $p,q,r$ are the same. The points
$(p,p,q)$ will not form a polar triangle if  ${}^t\!p A p\not= 0$ 
or ${}^t\!p A q \not= 0$.

\medskip
For the case ${}^t\!p A p\not= 0$  we pick coordinates such that $A$ is the 
identity matrix and $p=(1\,\colon 0\,\colon 0)$. Let
\[ B=\left\{ \begin{array}{ll}
\left({\textstyle \begin{array}{ccc}
0 & \!\!\!-q_2 &  q_1\\
\!\!\!-q_2 & 0 & 0\\
q_1 & 0 & 0
\end{array}}\right)
& \mbox{for\ } q_0\not=0 \mbox{\ or\ }q_1^2+q_2^2\not=0\\[4.5ex]
 \left({\textstyle \begin{array}{ccc}
0 & 0 &  1\\
0 & 1 & \!\!\pm i\\
1 & \!\!\pm i & \!-1
\end{array}}\right)
& \mbox{for\ } q=(0\,\colon  1\,\colon \pm i),
\end{array} \right. \]
then $C_B$ is a conic of ${\cal B}$ that contains $p$ and $q$,
but is smooth in $p$, and its 
tangent in $p$ is not the polar of $q$, so it does not contain $p$ twice.

\medskip
Finally, if ${}^t\! p A p=0$ and ${}^t\! p A q \not=0$, 
we pick coordinates such that $A$ is the 
identity matrix and $p=(1\,\colon i\,\colon 0)$. Let 
\[ B=\left\{ \begin{array}{ll}
 \left({\textstyle \begin{array}{ccc}
-2i & 2 & \!\!\! iq_0-2q_1\\
2 & 2i & -iq_1\\
\!\! iq_0-2q_1 &\!\! -iq_1 & 0
\end{array}}\right)
& \mbox{for\ } q=(q_0\,\colon  q_1\,\colon 1) \\[4.5ex]
\left({\textstyle \begin{array}{ccc}
0 & 0 &  1\\
0 & 0 & 0\\
1 & 0 & 0
\end{array}}\right)
& \mbox{for\ } q_2=0, 
\end{array} \right. \]
then we are in the same situation as above. 
\hfill $\Box$

\textsc{Mathematisches Institut IV, 
Heinrich-Heine-Universit\"at, 
Universit\"atsstr. 1,
40225 D\"usseldorf,
Germany}\\
\texttt{piontkow@uni-duesseldorf.de}

\medskip
\textsc{Faculteit der  Wiskunde en Natuurwetenschappen,
Afdeling Wiskunde en Informatica,
Postbus 9215,
2300 RA Leiden,
Netherlands}\\
\texttt{VEN@rulwinw.leidenuniv.nl}

\end{document}